\documentclass[11pt]{article}
\usepackage{mathrsfs}
\usepackage{amsfonts}
\usepackage{mathrsfs}
\usepackage{amsmath}
\usepackage{amssymb}
\usepackage{graphicx}
\graphicspath{{figure/}}
\usepackage{epic}
\usepackage{amsthm,latexsym,CJK}
\renewcommand{\paragraph}{\roman{paragraph}}
\def \n{\nabla}
\newtheorem{theorem}{Theorem}[section]
\newtheorem{remark}{Remark}[section]
\newtheorem{lemma}[theorem]{Lemma}
\newtheorem{corollary}[theorem]{Corollary}

\begin{document}
\title{\bf Matrix Li-Yau-Hamilton Estimates for Nonlinear Heat Equations\thanks{ Email: renx@cumt.edu.cn}}
\author{Xin-An Ren\\
\small Department of Mathematics, China University of Mining
and Technology, Xuzhou 221116, China.}

\date{}
\maketitle

\renewcommand{\theequation}{\thesection.\arabic{equation}}

\tableofcontents

\begin{abstract}
In this paper we are concerned with the matrix Li-Yau-Hamilton estimates for nonlinear heat equations. Firstly, we derive such estimate on a K\"{a}hler manifold with a fixed K\"{a}hler metric. Then we consider the estimate on K\"{a}hler manifolds with K\"{a}hler metrics evolving under the rescaled K\"{a}hler-Ricci flow. Both of the estimates are generalized to constrained cases. Finally, we extend the estimtes to more general nonlinear heat equations on both Riemannian manifolds and K\"{a}hler manifolds.
\end{abstract}

\section{Introduction}
In the seminar paper \cite{LY}, Li and Yau first developed the fundamental gradient estimates for any positive solutions of heat equation
 \begin{equation}
\frac{\partial}{\partial t}u=\Delta u
 \end{equation}
 on a Riemannian manifold, from which they derived the classical Harnack inequalities. Later, similar technique was employed to consider the geometric evolution equations. In this direction, Hamilton proved Harnack inequalities for Ricci flow in \cite{H3,H5} and mean curvature flow in \cite{H7}. So this type of Harnack inequalities was called Li-Yau-Hamilton estimate in \cite{NT1}. Similar Li-Yau-Hamilton estimates for Gauss curvature flow and Yamabe flow were proved by Chow in \cite{Ch1,Ch3}. Meanwhile, the estimate for K\"{a}hler-Ricci flow was proved by Cao in \cite{C2}. Later in \cite{A}, the corresponding estimate was established for a general class of hypersurface flows by Andrews. Moreover, a nice geometric interpretation of Hamilton quantities in the Li-Yau-Hamilton estimate for Ricci flow was given by Chow and Chu in \cite{CC1}. Using this geometric approach, Chow and Knopf gave a new Li-Yau-Hamilton estimate for Ricci flow in \cite{CK}.

  In another direction, Hamilton extended Li-Yau's gradient estimate of heat equation to the full matrix version in \cite{H4}. Precisely, he proved that a positive solution of heat equation
on a closed Riemannian manifold with nonnegative sectional curvature and parallel Ricci curvature satisfies the estimate
\begin{equation}
\nabla_i\nabla_j\ln u+\frac{1}{2t}g_{ij}\geq 0.
\end{equation}
Moreover, he derived several important monotonicity formulae out of the estimate in \cite{H6}.
 The assumption that $(M,g)$ has parallel Ricci curvature in Hamilton's estimate is rather restrictive and can be removed if the manifold is K\"{a}hler. This was observed by Cao and Ni for the K\"{a}hler manifold with a fixed metric in \cite{CN} and by Chow and Ni for the K\"{a}hler manifolds with time dependent metrics evolving under the K\"{a}hler-Ricci flow in \cite{N}. These two estimates can be unified by the interpolation method, which was originally proposed by Chow in \cite{Ch4}.
\vspace{2mm}\newline
{\bf Theorem A}(Chow-Ni) {\it Assume that $(M,g(t))$ satisfies $\epsilon$-K\"{a}hler-Ricci flow
 \begin{equation}
 \frac{\partial}{\partial t}g_{i\bar{j}}=-\epsilon R_{i\bar{j}}
 \end{equation}
 with nonnegative holomorphic bisectional curvature.
 In the case that $M$ is complete noncompact, additionally assume that the bisectional curvature is bounded. If $u$ is a positive solution to
\begin{equation}
\frac{\partial}{\partial t}u=\Delta u+\epsilon Ru,
\end{equation}
where $R$ is the scalar curvature, then
\begin{equation}
\nabla_i\nabla_{\bar{j}}\ln u+\epsilon R_{i\bar{j}}+\frac{1}{t}g_{i\bar{j}}\geq 0.
\end{equation}}

It should be pointed out that there are also other important works along this direction, see \cite{CC2,N2,N4,NT1,NT2}.

To better understand how the Li-Yau-Hamilton estimate for the Ricci flow can be perturbed or extended, Chow and Hamilton extended the matrix Li-Yau-Hamilton estimate on Riemannian manifolds to the constrained case in \cite{CH}.
\vspace{2mm}\newline
{\bf Theorem B}(Chow-Hamilton) {\it Let $(M,g)$ be a compact Riemannian manifold with nonnegative sectional curvature and parallel Ricci curvature. If $u$ and $v$ are two solutions of heat equation
\begin{equation}
\frac{\partial}{\partial t}u=\Delta u,\hspace{4mm}\frac{\partial}{\partial t}v=\Delta v,
\end{equation}
with $|v|<u$, then
\begin{equation}
\nabla_i\nabla_j\ln u+\frac{1}{2t}g_{ij}> \frac{\n_ih\n_{j}h}{1-h^2},
\end{equation}
where $h=v/u$.}\vspace{2mm}

In \cite{RYSZ}, we extended the matrix Li-Yau-Hamilton estimates due to Cao-Ni and Chow-Ni on K\"{a}hler manifolds to the constrained case.
\vspace{2mm}\newline
{\bf Theorem C}(Ren-Yao-Shen-Zhang) {\it
Let $(M,g(t))$ be a compact solution of $\epsilon$-K\"{a}hler-Ricci flow
\begin{equation}
 \frac{\partial}{\partial t}g_{i\bar{j}}=-\epsilon R_{i\bar{j}}
 \end{equation}
 with nonnegative holomorphic bisectional curvature.
 If $u$ and $v$ are two solutions of the equation
\begin{equation}
\frac{\partial}{\partial t}u=\Delta u+\epsilon Ru, \hspace{6mm}\frac{\partial}{\partial t}v=\Delta v+\epsilon Rv \label{heatequ}
\end{equation}
with $|v|<u$, then we have
\begin{equation}
\n_i\n_{\bar{j}}\ln u+\frac{1}{t}g_{i\bar{j}}+\epsilon R_{i\bar{j}}> \frac{\n_ih\n_{\bar{j}}h}{1-h^2},
\end{equation}
where $h=v/u$.}

In \cite{CFL}, the authors derived the Li-Yau-Hamilton estimate for nonlinear heat equation
\begin{equation}
\frac{\partial}{\partial t}u=\Delta u+au\ln u.
\end{equation}
Furthermore, the matrix version of this estimate was generalized by Wu in \cite{W}.
\vspace{2mm}\newline
{\bf Theorem D}(Wu) {\it
Let $(M,g)$ be a compact Riemannian manifold with nonnegative sectional curvature and parallel Ricci curvature. If $u$ is a positive solution to the equation
\begin{equation}
\frac{\partial}{\partial t}u=\Delta u+au\ln u,
\end{equation}
 then
\begin{equation}
\nabla_i\nabla_{j}\ln u+\frac{a}{2(1-e^{-at})}g_{ij}>0.
\end{equation}
}

In this paper, we firstly consider the matrix Li-Yau-Hamilton estimate for the nonlinear heat equation on a K\"{a}hler manifold and obtain the following
\begin{theorem}
Let $(M,g)$ be a compact K\"{a}hler manifold with nonnegative holomorphic bisectional curvature. If $u$ is a positive solution to the nonlinear heat equation
\begin{equation}
\frac{\partial}{\partial t}u=\Delta u+au\ln u, \label{NL}
\end{equation}
 then we have
\begin{equation}
\nabla_i\nabla_{\bar{j}}\ln u+\frac{a}{1-e^{-at}}g_{i\bar{j}}>0.
\end{equation}
\end{theorem}

In \cite{CK}, the Ricci flow rescaled by a cosmological constant was introduced by Chow and Knopf to prove a new Li-Yau-Hamilton estimate for the Ricci flow. Then Kotschwar used the similar technique to introduce the mean curvature flow with a cosmological constant in his thesis \cite{K}. In this paper the following evolution equation
\begin{equation}
\frac{\partial}{\partial t}g_{i\bar{j}}=-R_{i\bar{j}}+ag_{i\bar{j}}
\end{equation}
is called the rescaled K\"{a}hler-Ricci flow for simplicity. It is obvious that if $a=0$ then this flow reduces to K\"{a}hler-Ricci flow.
\begin{theorem}
Let $g(t)$ be a solution to the rescaled K\"{a}hler-Ricci flow
 \begin{equation}
\frac{\partial}{\partial t}g_{i\bar{j}}=-R_{i\bar{j}}+ag_{i\bar{j}}
\end{equation}
 on compact manifold $M$ with nonnegative holomorphic bisectional curvature. If $u$ is a positive solution to the nonlinear heat equation
\begin{equation}
\frac{\partial}{\partial t}u=\Delta u+ Ru+au\ln u, \label{NLF}
\end{equation}
where $a$ is a positive constant,
then we have
\begin{equation}
\nabla_i\nabla_{\bar{j}}\ln u+R_{i\bar{j}}+\frac{a}{1-e^{-at}}g_{i\bar{j}}>0.
\end{equation}
\end{theorem}

\begin{remark} By taking limits $a\rightarrow 0$ and $a\rightarrow 0^+$ in Theorem 1.1 and Theorem 1.2, respectively, we have that
 $$\frac{a}{1-e^{-at}}\rightarrow\frac{1}{t}.$$
 So Theorem 1.1 and Theorem 1.2 can be regarded as generalizations of theorems due to Cao-Ni in \cite{CN} and Chow-Ni in \cite{N}, respectively.
\end{remark}

Moreover, these results can also be extended to constrained cases.

\begin{theorem}
Let $(M,g)$ be a compact K\"{a}hler manifold. If $u$ and $v$ are positive solutions to the nonlinear heat equation
\begin{equation}
\frac{\partial}{\partial t}u=\Delta u+au\ln u
\end{equation}
and $h=v/u$, then in any of the two cases:
\newline (i) $a>0$, $0<h<1$ and holomorphic bisectional curvature is nonnegative,\\
(ii) $a<0$, $0<c<h<1$, where $c$ is a free parameter, and holomorphic bisectional curvature satisfies $R_{i\bar{j}k\bar{l}}\geq -aK(g_{i\bar{j}}g_{k\bar{l}}+g_{i\bar{l}}g_{k\bar{j}})$ with $K$ satisfying $K\geq -1-\frac{2\ln c}{1-c^2}$,\\
the following estimate holds
\begin{equation}
\nabla_i\nabla_{\bar{j}}\ln u+\frac{a}{1-e^{-at}}g_{i\bar{j}}>\frac{\n_ih\n_{\bar{j}}h}{1-h^2}.
\end{equation}
\end{theorem}

\begin{theorem}
Let $(M,g(t))$ be a compact solution to the rescaled K\"{a}hler-Ricci flow with nonnegative holomorphic bisectional curvature. If $u$ and $v$ are positive solutions to the nonlinear heat equation
\begin{equation}
\frac{\partial}{\partial t}u=\Delta u+ Ru+au\ln u
\end{equation}  with $v<u$ and $a>0$, then we have
\begin{equation}
\nabla_i\nabla_{\bar{j}}\ln u+ R_{i\bar{j}}+\frac{a}{1-e^{-at}}g_{i\bar{j}}>\frac{\n_ih\n_{\bar{j}}h}{1-h^2},
\end{equation}
where $h=v/u$.
\end{theorem}

In the proof of the theorems above, we let
\begin{equation}
L=\ln u.
\end{equation}
Then (\ref{NL}) turns out to be
\begin{equation}
\frac{\partial}{\partial t}L=\Delta L+|\n L|^2+aL.
\end{equation}
So it is interesting to turn our attention to more general equations.
\begin{theorem}
Let $(M,g)$ be a compact Riemannian manifold with nonnegative sectional curvature and parallel Ricci curvature and let $L$ be a solution to the nonlinear heat equation
\begin{equation}
\frac{\partial}{\partial t}L=\Delta L+|\nabla L|^2+F(L),
\end{equation}
where $F$ is a convex function of $L$, i.e., $F''(L)\geq0$. If
\begin{equation}
\nabla_i\nabla_j L+f(t)g_{ij}\geq 0
\end{equation}
at $t=0$, where $f(t)$ satisfies
\begin{equation}
2f^2(t)-F'(L)f(t)+f'(t)\geq 0,\label{fconditionRiemann}
\end{equation}
then it persists for $t>0$.
\end{theorem}

\begin{theorem}
Let $(M,g)$ be a compact K\"{a}hler manifold with nonnegative holomorphic bisectional curvature and let $L$ be a solution to the nonlinear heat equation
\begin{equation}
\frac{\partial}{\partial t}L=\Delta L+|\nabla L|^2+F(L),
\end{equation}
where $F$ is a convex function of $L$, i.e., $F''(L)\geq 0$. If
\begin{equation}
\nabla_i\nabla_{\bar j} L+f(t)g_{i{\bar j}}\geq 0,
\end{equation}
$t=0$, where $f(t)$ satisfies
\begin{equation}
f^2(t)-F'(L)f(t)+f'(t)\geq 0,\label{fconditionKahler}
\end{equation}
then it persists for $t>0$.
\end{theorem}

\begin{remark}
If $F(L)=aL$, then by solving the differential equation
$$2f^2(t)-af(t)+f'(t)=0,$$
we have that
$$f(t)=\frac{a}{2(1-e^{-at})}.$$
So Theorem 1.5 can be considered as a generalization of Theorem D. Similarly, Theorem 1.6 generalizes Theorem 1.1 when $F(L)=aL$.

\end{remark}

The rest of the paper is organized as follows. We devote section 2 to the notations and basic formulae in Riemannian geometry and K\"{a}hler geometry. In section 3, we prove Theorem 1.1 and Theorem by using Hamilton's maximum principle for tensors. In section 4, we consider the constrained case and prove Theorem 1.3 and Theorem 1.4. Finally, we prove Theorem 1.5 and Theorem 1.6 in section 5.

\section{Preliminaries}\setcounter{equation}{0}
Let $(M,g)$ be a Riemannnian manifold and $\nabla$ be the Levi-Civita connection of the metric $g$. Then Riemann curvature tensor is defined by
\begin{equation}
R(X,Y)Z=\nabla_X\nabla_YZ-\nabla_Y\nabla_X Z-\nabla_{[X,Y]}Z.
\end{equation}
In local coordinate, its components are determined by
\begin{equation}
R(\frac{\partial}{\partial x^i},\frac{\partial}{\partial x^j})\frac{\partial}{\partial x^k}=R^l_{ijk}\frac{\partial}{\partial x^l}.
\end{equation}
We keep the notation
\begin{equation}
R_{ijkl}=g_{lm}R^m_{ijk}
\end{equation}
as in \cite{CLN}.
The Ricci tensor is defined by
\begin{equation}
R_{ij}=R^k_{kij}
\end{equation}
and scalar curvature is given by
\begin{equation}
R=g^{ij}R_{ij}.
\end{equation}
We say that $M$ has nonnegative sectional curvature if
\begin{equation}
R_{kijl}v^iv^jw^kw^l\geq 0
\end{equation}
for all $v,w\in T_xM$ at $x\in M$.
We define Ricci curvature to be parallel if $\nabla_k R_{ij}=0$ for all $i,j,k$. To communicate covariant differentiation, we need following Ricci identities. If $\alpha$ is a $1$-form, then
\begin{equation}
\nabla_i\nabla_j\alpha_k=\nabla_j\nabla_i\alpha_k-R^l_{ijk}\alpha_l.
\end{equation}

Let $M$ be a  K\"{a}hler manifold with K\"{a}hler metric $g_{i\bar{j}}$. The K\"{a}hler form
\begin{equation}
\omega=\frac{\sqrt{-1}}{2}g_{i\bar{j}}dz^i\wedge d\bar{z}^j
\end{equation}
is a closed real $(1,1)$-form, so we have
\begin{equation}
\frac{\partial g_{i\bar{j}}}{\partial z^k}=\frac{\partial g_{k\bar{j}}}{\partial z^i},\hspace{4mm}\frac{\partial g_{i\bar{j}}}{\partial \bar{z}^k}=\frac{\partial g_{i\bar{k}}}{\partial \bar{z}^j}.
\end{equation}

The Christoffel symbols of the metric $g_{i\bar{j}}$ are given by
\begin{equation}
\Gamma^k_{ij}=g^{k\bar{l}}\frac{\partial g_{i\bar{l}}}{\partial z^j}, \hspace{4mm} \Gamma^{\bar{l}}_{\bar{i}\bar{j}}=g^{k\bar{l}}\frac{\partial g_{k\bar{j}}}{\partial \bar{z}^i},
\end{equation}
where $g^{i\bar{j}}=(g_{i\bar{j}})^{-1}$. It is easy to see that $\Gamma^k_{ij}$ is symmetric in $i$ and $j$ and $\Gamma^{\bar{l}}_{\bar{i}\bar{j}}$ is symmetric in $\bar{i}$ and $\bar{j}$.

The curvature tensor of the metric $g_{i\bar{j}}$ is defined as
\begin{equation}
R^j_{ik\bar{l}}=\frac{\partial\Gamma^j_{ik}}{\partial \bar{z}^l},\hspace{4mm}
R_{i\bar{j}k\bar{l}}=g_{p\bar{j}}R^p_{ik\bar{l}}.
\end{equation}
It is easy to see that $R_{i\bar{j}k\bar{l}}$ is symmetric in $i$ and $k$, in $\bar{j}$ and $\bar{l}$ and in pairs $i\bar{j}$ and $k\bar{l}$. The second Bianchi identity in K\"{a}hler case reduces to
\begin{equation}
\nabla_p R_{i\bar{j}k\bar{l}}=\nabla_k R_{i\bar{j}p\bar{l}}, \hspace{4mm}\nabla_{\bar{q}}R_{i\bar{j}k\bar{l}}=\nabla_{\bar{l}}R_{i\bar{j}k\bar{q}}.
\end{equation}

$M$ is said to have nonnegative holomorphic bisectional curvature if
\begin{equation}
R_{i\bar{j}k\bar{l}}v^iv^{\bar{j}}w^kw^{\bar{l}}\geq 0
\end{equation}
for all nonzero vectors $v$ and $w$ in the holomorphic tangent space $T_xM$ at $x\in M$.

The Ricci tensor of the metric $g_{i\bar{j}}$ is obtained by taking the trace of the curvature tensor:
\begin{equation}
R_{i\bar{j}}=g^{k\bar{l}}R_{{i\bar{j}k\bar{l}}}
\end{equation}
and the scalar curvature is given by
\begin{equation}
R=g^{i\bar{j}}R_{i\bar{j}}.
\end{equation}

Finally, we give the commutation formulae for covariant differentiations in K\"{a}hler geometry. Covariant differentiations of the same type can be commuted freely, e.g.,
\begin{equation}
\nabla_k\nabla_jv_i=\nabla_j\nabla_kv_i,\hspace{4mm}\nabla_{\bar{k}}\nabla_{\bar{j}}v_i=\nabla_{\bar{j}}\nabla_{\bar{k}}v_i.
\end{equation}
But we shall need following formulas when commuting covariant derivatives of different types
\begin{equation}
\nabla_k\nabla_{\bar{j}}v_i=\nabla_{\bar{j}}\nabla_kv_i-R_{k{\bar{j}}i{\bar{l}}}v_l.
\end{equation}
\section{Matrix Li-Yau-Hamilton estimates}\setcounter{equation}{0}
In this section, we will give the proof of Theorem 1.1 and Theorem 1.2. Firstly, we give some lemmas used in the proof. It should be pointed out that all the computations are taken in normal coordinate.

\begin{lemma}
If $A$ is a smooth function satisfying the evolution equation
\begin{equation}
\frac{\partial}{\partial t}A=\Delta {A}+B
\end{equation}
on a K\"{a}hler manifold, then $\nabla_i\nabla_{\bar{j}}A$ satisfies
\begin{equation}
\begin{split}
\frac{\partial}{\partial t}\left(\nabla_i\nabla_{\bar{j}}A\right)=&\Delta\left(\nabla_i\nabla_{\bar{j}}A\right)+R_{i\bar{j}l\bar{k}}\nabla_k\nabla_{\bar{l}}A\\&-\frac{1}{2}\left(R_{l\bar{j}}\nabla_i\nabla_{\bar{l}}A+R_{i\bar{l}}\nabla_l\nabla_{\bar{j}}A\right)
+\nabla_i\nabla_{\bar{j}}B.
\end{split}
\end{equation}
\end{lemma}
Proof: A straightforward computation leads to
\begin{equation}
\begin{split}
\frac{\partial}{\partial t}\left(\nabla_i\nabla_{\bar{j}}A\right)&=\nabla_i\nabla_{\bar{j}}\left(\frac{\partial}{\partial t}A\right)\\
&=\nabla_i\nabla_{\bar{j}}\left(\Delta {A}+B\right)\\
&=\nabla_i\nabla_{\bar{j}}\left(\Delta {A}\right)+\nabla_i\nabla_{\bar{j}}B.
\end{split}
\end{equation}
So it is sufficient to show that
\begin{equation}\label{evol1}
\nabla_i\nabla_{\bar{j}}\left(\Delta {A}\right)=\Delta\left(\nabla_i\nabla_{\bar{j}}A\right)+R_{i\bar{j}l\bar{k}}\nabla_k\nabla_{\bar{l}}A-\frac{1}{2}\left(R_{l\bar{j}}\nabla_i\nabla_{\bar{l}}A+R_{i\bar{l}}\nabla_l\nabla_{\bar{j}}A\right)
\end{equation}
Notice that
 $$\Delta {A}=\frac{1}{2}\left(\nabla_k\nabla_{\bar{k}}A+\nabla_{\bar{k}}\nabla_{k}A\right).$$
 By using communication formula, we get
\begin{equation}
\begin{split}
\nabla_i\nabla_{\bar{j}}\left(\nabla_k\nabla_{\bar{k}}A\right)=&\nabla_i\left(\nabla_k\nabla_{\bar{j}}\nabla_{\bar{k}}A-R_{l\bar{j}}\nabla_{\bar{l}}A\right)\\
=&\nabla_k\nabla_i\nabla_{\bar{k}}\nabla_{\bar{j}}A-\nabla_i{R_{l\bar{j}}}\nabla_{\bar{l}}A-R_{l\bar{j}}\nabla_i\nabla_{\bar{l}}A\\
=&\nabla_k\left(\nabla_{\bar{k}}\nabla_i\nabla_{\bar{j}}A+R_{i\bar{k}l\bar{j}}\nabla_{\bar{l}}A\right)-\nabla_i{R_{l\bar{j}}}\nabla_{\bar{l}}A-R_{l\bar{j}}\nabla_i\nabla_{\bar{l}}A\\
=&\nabla_k\nabla_{\bar{k}}\nabla_i\nabla_{\bar{j}}A+R_{i\bar{k}l\bar{j}}\nabla_k\nabla_{\bar{l}}A+\nabla_k{R_{i\bar{k}l\bar{j}}}\nabla_{\bar{l}}A
-\nabla_i{R_{l\bar{j}}}\nabla_{\bar{l}}A-R_{l\bar{j}}\nabla_i\nabla_{\bar{l}}A\\
=&\nabla_k\nabla_{\bar{k}}\nabla_i\nabla_{\bar{j}}A+R_{i\bar{j}l\bar{k}}\nabla_k\nabla_{\bar{l}}A-R_{l\bar{j}}\nabla_i\nabla_{\bar{l}}A  \label{evol2}
\end{split}
\end{equation}
and
\begin{equation}
\begin{split}
\nabla_i\nabla_{\bar{j}}\left(\nabla_{\bar{k}}\nabla_{k}A\right)=&\nabla_i\left(\nabla_{\bar{k}}\nabla_{\bar{j}}\nabla_{k}A\right)\\=&\nabla_i\left(\nabla_{\bar{k}}\nabla_{k}\nabla_{\bar{j}}A\right)\\
=&\nabla_{\bar{k}}\nabla_i\nabla_k\nabla_{\bar{j}}A+R_{i\bar{k}l\bar{j}}\nabla_k\nabla_{\bar{l}}A-R_{i\bar{k}k\bar{l}}\nabla_l\nabla_{\bar{j}}A\\
=&\nabla_{\bar{k}}\nabla_k\left(\nabla_i\nabla_{\bar{j}}A\right)+R_{i\bar{j}l\bar{k}}\nabla_k\nabla_{\bar{l}}A-R_{i\bar{l}}\nabla_l\nabla_{\bar{j}}A. \label{evol3}
\end{split}
\end{equation}
Inserting (\ref{evol2}) and (\ref{evol3}) into (\ref{evol1}), we obtain the evolution equation of $\nabla_i\nabla_{\bar{j}}A$.

Before proceeding further, let us remark that the evolution equation (3.2) is still valid if the K\"{a}hler metric is moving under the rescaled K\"{a}hler-Ricci flow.

\begin{lemma} Assume that $u$ is a positive solution to the nonlinear heat equation
 $$\frac{\partial}{\partial t}u=\Delta u+au\ln u.$$
 If we set $L=\ln u$, then the evolution equation of $\n_{i}\n_{\bar{j}}L$ is
\begin{equation}
\begin{split}
\frac{\partial}{\partial t}\nabla_i\nabla_{\bar{j}}L=&\Delta\left(\nabla_i\nabla_{\bar{j}}L\right)+R_{i\bar{j}l\bar{k}}\nabla_k\nabla_{\bar{l}}L
+R_{i\bar{j}k\bar{l}}{\nabla_{l}L\nabla_{\bar{k}}L}\\&+\nabla_{k}L\nabla_{\bar{k}}\nabla_i\nabla_{\bar{j}}L+\nabla_{\bar{k}}L\nabla_k\nabla_i\nabla_{\bar{j}}L
\\&+\nabla_{\bar{j}}\nabla_{k}L\nabla_i\nabla_{\bar{k}}L+\nabla_i\nabla_{k}L\nabla_{\bar{j}}\nabla_{\bar{k}}L
\\&-\frac{1}{2}\left(R_{l\bar{j}}\nabla_i\nabla_{\bar{l}}L+R_{i\bar{l}}\nabla_l\nabla_{\bar{j}}L\right)+a\nabla_i\nabla_{\bar{j}}L.
\end{split}
\end{equation}
\end{lemma}
Proof: We use
\begin{equation}
\frac{\partial}{\partial t}L=\frac{\partial}{\partial t}\ln u=\frac{u_t}{u}=\frac{1}{u}\left(\Delta u+au\ln u\right)=\frac{\Delta u}{u}+a\ln u
\end{equation}
and
\begin{equation}
\Delta L=\n_i\n_{\bar{i}}L=\n_i\frac{\n_{\bar{i}}u}{u}=\frac{\Delta u}{u}-\frac{|\nabla u|^2}{u^2}
\end{equation}
to conclude that
\begin{equation}
\frac{\partial}{\partial t}L=\Delta L+\frac{|\nabla u|^2}{u^2}+aL=\Delta L+|\n L|^2+aL.
\end{equation}

Applying Lemma 3.1 to $A=L$, we obtain
\begin{equation}
\begin{split}
\frac{\partial}{\partial t}\left(\nabla_i\nabla_{\bar{j}}L\right)=&\Delta\left(\nabla_i\nabla_{\bar{j}}L\right)+R_{i\bar{j}l\bar{k}}\nabla_k\nabla_{\bar{l}}L-\frac{1}{2}\left(R_{l\bar{j}}\nabla_i\nabla_{\bar{l}}L+R_{i\bar{l}}\nabla_l\nabla_{\bar{j}}L\right)\\
&+\nabla_i\nabla_{\bar{j}}|\nabla L|^2
+a\nabla_i\nabla_{\bar{j}}L.
\end{split}
\end{equation}
Now the lemma follows from the following computation
\begin{equation}
\begin{split}
\nabla_i\nabla_{\bar{j}}|\nabla L|^2=&\nabla_i\nabla_{\bar{j}}\left(\nabla_{k}L\nabla_{\bar{k}}L\right)\\=&\nabla_i\left(\nabla_{\bar{j}}\nabla_{k}L\nabla_{\bar{k}}L+\nabla_{k}L\nabla_{\bar{j}}\nabla_{\bar{k}}L\right)
\\=&\nabla_i\nabla_{\bar{j}}\nabla_{k}L\nabla_{\bar{k}}L+\nabla_{\bar{j}}\nabla_{k}L\nabla_i\nabla_{\bar{k}}L+\nabla_i\nabla_{k}L\nabla_{\bar{j}}\nabla_{\bar{k}}L
+\nabla_{k}L\nabla_i\nabla_{\bar{j}}\nabla_{\bar{k}}L
\\=&\nabla_{\bar{k}}L\nabla_k\left(\nabla_i\nabla_{\bar{j}}L\right)+\nabla_{\bar{j}}\nabla_{k}L\nabla_i\nabla_{\bar{k}}L\\&+\nabla_i\nabla_{k}L\nabla_{\bar{j}}\nabla_{\bar{k}}L
+\nabla_{k}L\nabla_{\bar{k}}\left(\nabla_i\nabla_{\bar{j}}L\right)+R_{i\bar{j}k\bar{l}}\nabla_{l}L\nabla_{\bar{k}}L.
\end{split}
\end{equation}

Now we are in a position to prove Theorem 1.1.

{\bf Proof of Theorem 1.1} Set
\begin{equation}
P_{i\bar{j}}=\nabla_i\nabla_{\bar{j}}L+\frac{a}{1-e^{-at}}g_{i\bar{j}}.
\end{equation}
We compute the evolution equation of $P_{i\bar{j}}$
\begin{equation}
\begin{split}
\frac{\partial}{\partial t}P_{i\bar{j}}=&\frac{\partial}{\partial t}\nabla_i\nabla_{\bar{j}}L-\frac{a^2e^{-at}}{(1-e^{-at})^2}g_{i\bar{j}}\\
=&\Delta\left(\nabla_i\nabla_{\bar{j}}L\right)+R_{i\bar{j}l\bar{k}}\nabla_k\nabla_{\bar{l}}L
+R_{i\bar{j}k\bar{l}}{\nabla_{l}L\nabla_{\bar{k}}L}\\&+\nabla_{k}L\nabla_{\bar{k}}\nabla_i\nabla_{\bar{j}}L+\nabla_{\bar{k}}L\nabla_k\nabla_i\nabla_{\bar{j}}L
\\&+\nabla_{\bar{j}}\nabla_{k}L\nabla_i\nabla_{\bar{k}}L+\nabla_i\nabla_{k}L\nabla_{\bar{j}}\nabla_{\bar{k}}L
\\&-\frac{1}{2}\left(R_{l\bar{j}}\nabla_i\nabla_{\bar{l}}L+R_{i\bar{l}}\nabla_l\nabla_{\bar{j}}L\right)+a\nabla_i\nabla_{\bar{j}}L-\frac{a^2e^{-at}}{(1-e^{-at})^2}g_{i\bar{j}}\\
=&\Delta P_{i\bar{j}}+R_{i\bar{j}l\bar{k}}P_{k\bar{l}}-\frac{a}{1-e^{-at}}R_{i\bar{j}}
+R_{i\bar{j}k\bar{l}}{\nabla_{l}L\nabla_{\bar{k}}L}+\nabla_{k}L\nabla_{\bar{k}}P_{i\bar{j}}+\nabla_{\bar{k}}L\nabla_kP_{i\bar{j}}\\
&+P_{i\bar{k}}P_{k\bar{j}}-\frac{2a}{1-e^{-at}}P_{i\bar{j}}+\frac{a^2}{(1-e^{-at})^2}g_{i\bar{j}}+\nabla_i\nabla_{k}L\nabla_{\bar{j}}\nabla_{\bar{k}}L\\
&-\frac{1}{2}\left(R_{l\bar{j}}P_{i\bar{l}}+R_{i\bar{l}}P_{l\bar{j}}\right)+\frac{a}{1-e^{-at}}R_{i\bar{j}}
+aP_{i\bar{j}}-\frac{a^2}{1-e^{-at}}g_{i\bar{j}}-\frac{a^2e^{-at}}{(1-e^{-at})^2}g_{i\bar{j}}\\
=&\Delta P_{i\bar{j}}+R_{i\bar{j}l\bar{k}}P_{k\bar{l}}
+R_{i\bar{j}k\bar{l}}{\nabla_{l}L\nabla_{\bar{k}}L}+\nabla_{k}L\nabla_{\bar{k}}P_{i\bar{j}}+\nabla_{\bar{k}}L\nabla_kP_{i\bar{j}}\\
&+P_{i\bar{k}}P_{k\bar{j}}-\frac{a(1+e^{-at})}{1-e^{-at}}P_{i\bar{j}}+\nabla_i\nabla_{k}L\nabla_{\bar{j}}\nabla_{\bar{k}}L-\frac{1}{2}\left(R_{l\bar{j}}P_{i\bar{l}}+R_{i\bar{l}}P_{l\bar{j}}\right)\\
\geq&\Delta P_{i\bar{j}}+\nabla_{k}L\nabla_{\bar{k}}P_{i\bar{j}}+\nabla_{\bar{k}}L\nabla_kP_{i\bar{j}}+R_{i\bar{j}l\bar{k}}P_{k\bar{l}}+P_{i\bar{k}}P_{k\bar{j}}-\frac{a(1+e^{-at})}{1-e^{-at}}P_{i\bar{j}}\\
&-\frac{1}{2}\left(R_{l\bar{j}}P_{i\bar{l}}+R_{i\bar{l}}P_{l\bar{j}}\right)\label{th11}
\end{split}
\end{equation}
In the last step, we have used the condition that the holomorphic bisectional curvature is nonnegative and the fact that
$\nabla_i\nabla_{k}L\nabla_{\bar{j}}\nabla_{\bar{k}}L$
is nonnegative definite. It is easy to see that
 $$R_{i\bar{j}l\bar{k}}P_{k\bar{l}}+P_{i\bar{k}}P_{k\bar{j}}-\frac{a(1+e^{-at})}{1-e^{-at}}P_{i\bar{j}}-\frac{1}{2}\left(R_{l\bar{j}}P_{i\bar{l}}+R_{i\bar{l}}P_{l\bar{j}}\right)$$
 satisfies the null-eigenvector condition of the tensor maximum principle. Hence we complete the proof of Theorem 1.1 by Hamilton's maximum principle for tensors.

 Now we turn to the proof of Theorem 1.2. The proof consists of the following lemmas.

\begin{lemma}(see \cite{C1}) If K\"{a}hler metrics $g_{i\bar{j}}(t)$ satisfy the rescaled K\"{a}hler-Ricci flow
\begin{equation}
\frac{\partial}{\partial t}g_{i\bar{j}}=-R_{i\bar{j}}+ag_{i\bar{j}},
\end{equation}
 then the Ricci tensors satisfy
\begin{equation}
\frac{\partial}{\partial t} R_{i\bar{j}}=\nabla_i\nabla_{\bar{j}}R=\Delta R_{i{\bar{j}}}+R_{i{\bar{j}}k{\bar{l}}}R_{l{\bar{k}}}-R_{i{\bar{p}}}R_{p{\bar{j}}}.
\end{equation}
\end{lemma}

 Following lemma is an easy consequence of Cao's Li-Yau-Hamilton esitmate for the K\"{a}hler-Ricci flow \cite{C2}.

\begin{lemma} If $(M,g(t))$ be a solution of the rescaled K\"{a}hler-Ricci flow
 \begin{equation}
\frac{\partial}{\partial t}g_{i\bar{j}}=-R_{i\bar{j}}+ag_{i\bar{j}}\label{NKRF}
\end{equation}
with nonnegative holomorphic bisectional curvature, then
\begin{equation}
\frac{\partial}{\partial t}R_{i\bar{j}}+R_{i\bar{k}}R_{k\bar{j}}+R_{i\bar{j},k}X^k+R_{i\bar{j},\bar{k}}X^{\bar{k}}
+R_{i\bar{j}k\bar{l}}X^kX^{\bar{l}}+\frac{ae^{-at}}{1-e^{-at}}R_{i\bar{j}}\geq 0.
\end{equation}
for any holomorphic vector $X$.
\end{lemma}
Proof:
Let $\hat{g}_{i\bar{j}}(s)$ be a solution to the K\"{a}hler-Ricci flow
\begin{equation}
\frac{\partial}{\partial s}\hat{g}_{i\bar{j}}(s)=-\hat{R}_{i\bar{j}}(s).\label{krf}
\end{equation}
Setting
\begin{equation}
s=\frac{1}{a}(1-e^{-at}), \hspace{3mm}g_{i\bar{j}}(t)=e^{at}\hat{g}_{i\bar{j}}(s),\label{mr1}
\end{equation}
we have that
\begin{equation}
\begin{split}
\frac{\partial}{\partial t}g_{i\bar{j}}(t)&=\frac{\partial}{\partial t}\left(e^{at}\hat{g}_{i\bar{j}}(s)\right)\\
&=ae^{at}\hat{g}_{i\bar{j}}(s)+e^{at}\frac{\partial}{\partial s}\hat{g}_{i\bar{j}}(s)\frac{\partial s}{\partial t}\\
&=ag_{i\bar{j}}(t)-R_{i\bar{j}}(t),
\end{split}
\end{equation}
which means that $g_{i\bar{j}}(t)=e^{at}\hat{g}_{i\bar{j}}(s)$ is a solution to the rescaled K\"{a}hler-Ricci flow (\ref{NKRF}).
The Li-Yau-Hamilton estimate for the  K\"{a}hler-Ricci flow (\ref{krf}) is
\begin{equation}
\frac{\partial }{\partial s}\hat{R}_{i\bar{j}}+\hat{R}_{i\bar{k}}\hat{R}_{k\bar{j}}+\hat{R}_{i\bar{j},k}v^k+\hat{R}_{i\bar{j},\bar{k}}v^{\bar{k}}
+\hat{R}_{i\bar{j}k\bar{l}}v^kv^{\bar{l}}+\frac{1}{s}\hat{R}_{i\bar{j}}\geq 0.
\end{equation}
If follows from (\ref{mr1}) that
\begin{equation}
\hat{g}^{i\bar{j}}=e^{at}g^{i\bar{j}},\hspace{2mm}\hat{R}_{i\bar{j}k\bar{l}}=e^{-at}R_{i\bar{j}k\bar{l}},\hspace{2mm}\hat{R}_{i\bar{j}}=R_{i\bar{j}}.
\end{equation}
Hence, we have
\begin{equation}
\frac{\partial}{\partial t}R_{i\bar{j}}\frac{\partial t}{\partial s}+e^{at}R_{i\bar{k}}R_{k\bar{j}}+R_{i\bar{j},k}v^k+R_{i\bar{j},\bar{k}}v^{\bar{k}}
+e^{-at}R_{i\bar{j}k\bar{l}}v^kv^{\bar{l}}+\frac{a}{1-e^{-at}}R_{i\bar{j}}\geq 0.\label{3.24}
\end{equation}
Multiplying $e^{-at}$ on the both side of (\ref{3.24}) and letting $X^k=e^{-at}v^k$, we get
\begin{equation}
\frac{\partial}{\partial t}R_{i\bar{j}}+R_{i\bar{k}}R_{k\bar{j}}+R_{i\bar{j},k}X^k+R_{i\bar{j},\bar{k}}X^{\bar{k}}
+R_{i\bar{j}k\bar{l}}X^kX^{\bar{l}}+\frac{ae^{-at}}{1-e^{-at}}R_{i\bar{j}}\geq 0.
\end{equation}
Now the proof of the lemma is completed.

\begin{lemma}Let $(M,g(t))$ be a solution to the rescaled K\"{a}hler-Ricci flow
 \begin{equation}
 \frac{\partial }{\partial t}g_{i\bar{j}}=-R_{i\bar{j}}+ag_{i\bar{j}}
 \end{equation}
 and let $u$ be a positive solution to the nonlinear heat equation
 $$\frac{\partial}{\partial t}u=\Delta u+Ru+au\ln u.$$
 If we set $L=\ln u$, then the evolution equation of $\n_{i}\n_{\bar{j}}L$ is
\begin{equation}
\begin{split}
\frac{\partial}{\partial t}\nabla_i\nabla_{\bar{j}}L=&\Delta\left(\nabla_i\nabla_{\bar{j}}L\right)+R_{i\bar{j}l\bar{k}}\nabla_k\nabla_{\bar{l}}L
+R_{i\bar{j}k\bar{l}}{\nabla_{l}L\nabla_{\bar{k}}L}\\&+\nabla_{k}L\nabla_{\bar{k}}\nabla_i\nabla_{\bar{j}}L+\nabla_{\bar{k}}L\nabla_k\nabla_i\nabla_{\bar{j}}L
\\&+\nabla_{\bar{j}}\nabla_{k}L\nabla_i\nabla_{\bar{k}}L+\nabla_i\nabla_{k}L\nabla_{\bar{j}}\nabla_{\bar{k}}L
\\&-\frac{1}{2}\left(R_{l\bar{j}}\nabla_i\nabla_{\bar{l}}L+R_{i\bar{l}}\nabla_l\nabla_{\bar{j}}L\right)+a\nabla_i\nabla_{\bar{j}}L\\
&+\Delta R_{i{\bar{j}}}+R_{i{\bar{j}}k{\bar{l}}}R_{l{\bar{k}}}-R_{i{\bar{p}}}R_{p{\bar{j}}}.
\end{split}
\end{equation}
\end{lemma}
Proof: By
\begin{equation}
\frac{\partial}{\partial t}L=\frac{\partial}{\partial t}\ln u=\frac{u_t}{u}=\frac{1}{u}\left(\Delta u+Ru+au\ln u\right)=\frac{\Delta u}{u}+R+a\ln u
\end{equation}
and
\begin{equation}
\Delta L=\n_i\n_{\bar{i}}L=\n_i\frac{\n_{\bar{i}}u}{u}=\frac{\Delta u}{u}-\frac{|\nabla u|^2}{u^2},
\end{equation}
we get the evolution equation of $L$
\begin{equation}
\frac{\partial}{\partial t}L=\Delta L+\frac{|\nabla u|^2}{u^2}+R+aL=\Delta L+|\n L|^2+R+aL.
\end{equation}
Applying Lemma 3.1 to $A=L$, we obtain
\begin{equation}
\begin{split}
\frac{\partial}{\partial t}\left(\nabla_i\nabla_{\bar{j}}L\right)=&\Delta\left(\nabla_i\nabla_{\bar{j}}L\right)+R_{i\bar{j}l\bar{k}}\nabla_k\nabla_{\bar{l}}L-\frac{1}{2}\left(R_{l\bar{j}}\nabla_i\nabla_{\bar{l}}L+R_{i\bar{l}}\nabla_l\nabla_{\bar{j}}L\right)\\
&+\nabla_i\nabla_{\bar{j}}|\nabla L|^2
+a\nabla_i\nabla_{\bar{j}}L+\nabla_i\nabla_{\bar{j}}R.
\end{split}
\end{equation}
Now the lemma follows from Lemma 3.3 and the following computation
\begin{equation}
\begin{split}
\nabla_i\nabla_{\bar{j}}|\nabla L|^2=&\nabla_i\nabla_{\bar{j}}\left(\nabla_{k}L\nabla_{\bar{k}}L\right)\\=&\nabla_i\left(\nabla_{\bar{j}}\nabla_{k}L\nabla_{\bar{k}}L+\nabla_{k}L\nabla_{\bar{j}}\nabla_{\bar{k}}L\right)
\\=&\nabla_i\nabla_{\bar{j}}\nabla_{k}L\nabla_{\bar{k}}L+\nabla_{\bar{j}}\nabla_{k}L\nabla_i\nabla_{\bar{k}}L+\nabla_i\nabla_{k}L\nabla_{\bar{j}}\nabla_{\bar{k}}L
+\nabla_{k}L\nabla_i\nabla_{\bar{j}}\nabla_{\bar{k}}L
\\=&\nabla_{\bar{k}}L\nabla_k\left(\nabla_i\nabla_{\bar{j}}L\right)+\nabla_{\bar{j}}\nabla_{k}L\nabla_i\nabla_{\bar{k}}L\\&+\nabla_i\nabla_{k}L\nabla_{\bar{j}}\nabla_{\bar{k}}L
+\nabla_{k}L\nabla_{\bar{k}}\left(\nabla_i\nabla_{\bar{j}}L\right)+R_{i\bar{j}k\bar{l}}\nabla_{l}L\nabla_{\bar{k}}L.
\end{split}
\end{equation}

Now we are in a position to prove Theorem 1.2

{\bf Proof of Theorem 1.2} Let
\begin{equation}
P_{i\bar{j}}=\nabla_{i}\nabla_{\bar{j}}L+R_{i\bar{j}}.
\end{equation}
It follows from Lemma 3.3 and Lemma 3.5 that
\begin{equation}
\begin{split}
\frac{\partial}{\partial t}P_{i\bar{j}}=&\frac{\partial}{\partial t}(\nabla_{i}\nabla_{\bar{j}}L+R_{i\bar{j}})\\
=&\Delta\left(\nabla_i\nabla_{\bar{j}}L\right)+R_{i\bar{j}l\bar{k}}\nabla_k\nabla_{\bar{l}}L
+R_{i\bar{j}k\bar{l}}{\nabla_{l}L\nabla_{\bar{k}}L}\\&+\nabla_{k}L\nabla_{\bar{k}}\nabla_i\nabla_{\bar{j}}L+\nabla_{\bar{k}}L\nabla_k\nabla_i\nabla_{\bar{j}}L
\\&+\nabla_{\bar{j}}\nabla_{k}L\nabla_i\nabla_{\bar{k}}L+\nabla_i\nabla_{k}L\nabla_{\bar{j}}\nabla_{\bar{k}}L
\\&-\frac{1}{2}\left(R_{l\bar{j}}\nabla_i\nabla_{\bar{l}}L+R_{i\bar{l}}\nabla_l\nabla_{\bar{j}}L\right)+a\nabla_i\nabla_{\bar{j}}L\\
&+2\Delta R_{i{\bar{j}}}+2R_{i{\bar{j}}k{\bar{l}}}R_{l{\bar{k}}}-2R_{i{\bar{p}}}R_{p{\bar{j}}}\\
=&\Delta P_{i{\bar{j}}}+R_{i\bar{j}l\bar{k}}P_{k{\bar{l}}}+\nabla_{k}L\nabla_{\bar{k}}P_{i{\bar{j}}}+\nabla_{\bar{k}}L\nabla_kP_{i{\bar{j}}}+(P_{i{\bar{k}}}-R_{i{\bar{k}}})(P_{k{\bar{j}}}-R_{k{\bar{j}}})\\
&+\nabla_i\nabla_{k}L\nabla_{\bar{j}}\nabla_{\bar{k}}L
-\frac{1}{2}\left(R_{l\bar{j}}P_{i{\bar{l}}}+R_{i\bar{l}}P_{l{\bar{j}}}\right)+aP_{i{\bar{j}}}-R_{i{\bar{p}}}R_{p{\bar{j}}}-aR_{i{\bar{j}}}\\
&+\Delta R_{i{\bar{j}}}+R_{i{\bar{j}}k{\bar{l}}}R_{l{\bar{k}}}-\nabla_{k}L\nabla_{\bar{k}}R_{i{\bar{j}}}-\nabla_{\bar{k}}L\nabla_kR_{i{\bar{j}}}+R_{i\bar{j}k\bar{l}}{\nabla_{l}L\nabla_{\bar{k}}L}\\
\geq&\Delta P_{i{\bar{j}}}+R_{i\bar{j}l\bar{k}}P_{k{\bar{l}}}+\nabla_{k}L\nabla_{\bar{k}}P_{i{\bar{j}}}+\nabla_{\bar{k}}L\nabla_kP_{i{\bar{j}}}+P_{i{\bar{k}}}P_{k{\bar{j}}}\\
&-\frac{3}{2}\left(R_{l\bar{j}}P_{i{\bar{l}}}+R_{i\bar{l}}P_{l{\bar{j}}}\right)+aP_{i{\bar{j}}}-aR_{i\bar{j}}-\frac{ae^{-at}}{1-e^{-at}}R_{i\bar{j}}.
\end{split}
\end{equation}
Set
\begin{equation}
\begin{split}
Q_{i{\bar{j}}}=P_{i{\bar{j}}}+\frac{a}{1-e^{-at}}g_{i{\bar{j}}}.
\end{split}
\end{equation}
Direct calculation gives rise to
\begin{equation}
\begin{split}
\frac{\partial}{\partial t}Q_{i{\bar{j}}}=&\frac{\partial}{\partial t}P_{i{\bar{j}}}+\frac{\partial}{\partial t}\left(\frac{a}{1-e^{-at}}\right)g_{i{\bar{j}}}+\frac{a}{1-e^{-at}}\frac{\partial}{\partial t}g_{i{\bar{j}}}\\
\geq&\Delta P_{i{\bar{j}}}+R_{i\bar{j}l\bar{k}}P_{k{\bar{l}}}+\nabla_{k}L\nabla_{\bar{k}}P_{i{\bar{j}}}+\nabla_{\bar{k}}L\nabla_kP_{i{\bar{j}}}+P_{i{\bar{k}}}P_{k{\bar{j}}}\\
&-\frac{3}{2}\left(R_{l\bar{j}}P_{i{\bar{l}}}+R_{i\bar{l}}P_{l{\bar{j}}}\right)+aP_{i{\bar{j}}}-\frac{a}{1-e^{-at}}R_{i\bar{j}}\\
&-\frac{a^2e^{-at}}{(1-e^{-at})^2}g_{i{\bar{j}}}+\frac{a}{1-e^{-at}}(-R_{i\bar{j}}+ag_{i\bar{j}})\\
=&\Delta Q_{i{\bar{j}}}+R_{i\bar{j}l\bar{k}}Q_{k{\bar{l}}}-\frac{a}{1-e^{-at}}R_{i{\bar{j}}}+\nabla_{k}L\nabla_{\bar{k}}Q_{i{\bar{j}}}+\nabla_{\bar{k}}L\nabla_kQ_{i{\bar{j}}}\\
&+(P_{i{\bar{k}}}+\frac{a}{1-e^{-at}}g_{i{\bar{k}}})(P_{k{\bar{j}}}-\frac{a}{1-e^{-at}}g_{k{\bar{j}}})+\frac{a^2}{(1-e^{-at})^2}g_{i{\bar{j}}}\\
&-\frac{3}{2}\left(R_{l\bar{j}}Q_{i{\bar{l}}}+R_{i\bar{l}}Q_{l{\bar{j}}}\right)+\frac{3a}{1-e^{-at}}R_{i{\bar{j}}}+aP_{i{\bar{j}}}-\frac{a}{1-e^{-at}}R_{i\bar{j}}\\
&-\frac{a^2e^{-at}}{(1-e^{-at})^2}g_{i{\bar{j}}}+\frac{a}{1-e^{-at}}(-R_{i\bar{j}}+ag_{i\bar{j}})\\
=&\Delta Q_{i{\bar{j}}}+R_{i\bar{j}l\bar{k}}Q_{k{\bar{l}}}+\nabla_{k}L\nabla_{\bar{k}}Q_{i{\bar{j}}}+\nabla_{\bar{k}}L\nabla_kQ_{i{\bar{j}}}\\
&+Q_{i{\bar{k}}}(P_{k{\bar{j}}}-\frac{a}{1-e^{-at}}g_{k{\bar{j}}})-\frac{3}{2}\left(R_{l\bar{j}}Q_{i{\bar{l}}}+R_{i\bar{l}}Q_{l{\bar{j}}}\right)+aQ_{i{\bar{j}}}+\frac{a^2}{1-e^{-at}}g_{i{\bar{j}}}\\
\geq&\Delta Q_{i{\bar{j}}}+\nabla_{k}L\nabla_{\bar{k}}Q_{i{\bar{j}}}+\nabla_{\bar{k}}L\nabla_kQ_{i{\bar{j}}}+R_{i\bar{j}l\bar{k}}Q_{k{\bar{l}}}\\
&+Q_{i{\bar{k}}}(P_{k{\bar{j}}}-\frac{a}{1-e^{-at}}g_{k{\bar{j}}})-\frac{3}{2}\left(R_{l\bar{j}}Q_{i{\bar{l}}}+R_{i\bar{l}}Q_{l{\bar{j}}}\right)+aQ_{i{\bar{j}}}
\end{split}
\end{equation}
It is obvious that
 $$R_{i\bar{j}l\bar{k}}Q_{k{\bar{l}}}+Q_{i{\bar{k}}}(P_{k{\bar{j}}}-\frac{a}{1-e^{-at}}g_{k{\bar{j}}})-\frac{3}{2}\left(R_{l\bar{j}}Q_{i{\bar{l}}}+R_{i\bar{l}}Q_{l{\bar{j}}}\right)+aQ_{i{\bar{j}}}
$$
 satisfies the null-eigenvector condition of the tensor maximum principle. Hence we complete the proof of Theorem 1.2.

 The following corollary, which follows easily from the relation (\ref{krf}) and (\ref{mr1}), provides an equivalent form of Theorem 1.2.
 \begin{corollary}
 Let $M$ be a compact K\"{a}hler manifold and $\hat{g}_{i\bar{j}}(s)$ be a solution to the K\"{a}hler-Ricci flow with nonnegative holomorphic bisectional curvature. If $u$ is a positive solution to the nonlinear heat equation (\ref{NLF}), then we have
 \begin{equation}
 \nabla_i\nabla_{\bar{j}}\ln u+\hat{R}_{i\bar{j}}+(\frac{1}{s}-a)\hat{g}_{i\bar{j}}>0,
 \end{equation}
 where $\hat{R}_{i\bar{j}}$ is the Ricci tensor of $\hat{g}_{i\bar{j}}(s)$.
 \end{corollary}
\section{Constrained matrix Li-Yau-Hamilton estimates}
\setcounter{equation}{0}
In this section, we will generalize Theorem 1.1 and Theorem 1.2 to constrained cases. The following lemma plays an important role in the proof of the theorems.
\begin{lemma}If $u$ and $v$ are positive solutions to the nonlinear heat equation (\ref{NL}), then the radio $h=v/u$ satisfies
\begin{equation}
\begin{split}
\frac{\partial}{\partial t}\left(\frac{\n_{i}h\n_{\bar{j}}h}{1-h^2}\right)
=&\Delta\left(\frac{\n_{i}h\n_{\bar{j}}h}{1-h^2}\right)+\n_{k}L\n_{\bar{k}}\left(\frac{\n_{i}h\n_{\bar{j}}h}{1-h^2}\right)+\n_{\bar{k}}L\n_{k}\left(\frac{\n_{i}h\n_{\bar{j}}h}{1-h^2}\right)\\
&-\frac{1}{1-h^2}\left(\n_{i}\n_{k}h+\frac{2h\n_{i}h\n_{k}h}{1-h^2}\right)\left(\n_{\bar{j}}\n_{\bar{k}}h+\frac{2h\n_{\bar{j}}h\n_{\bar{k}}h}{1-h^2}\right)\\
&-\frac{1}{1-h^2}\left(\n_{i}\n_{\bar{k}}h+\frac{2h\n_{i}h\n_{\bar{k}}h}{1-h^2}\right)\left(\n_{\bar{j}}\n_{k}h+\frac{2h\n_{\bar{j}}h\n_{k}h}{1-h^2}\right)\\
&+\n_{i}\n_{k}L\left(\frac{\n_{\bar{j}}h\n_{\bar{k}}h}{1-h^2}\right)+\n_{i}\n_{\bar{k}}L\left(\frac{\n_{\bar{j}}h\n_{k}h}{1-h^2}\right)\\&+\n_{\bar{j}}\n_{k}L\left(\frac{\n_{i}h\n_{\bar{k}}h}{1-h^2}\right)
+\n_{\bar{j}}\n_{\bar{k}}L\left(\frac{\n_{i}h\n_{k}h}{1-h^2}\right)
\\&-\frac{1}{2}R_{i\bar{l}}\left(\frac{\n_{\bar{j}}h\n_{l}h}{1-h^2}\right)-\frac{1}{2}R_{k\bar{j}}\left(\frac{\n_{\bar{k}}h\n_{i}h}{1-h^2}\right)\\&-\frac{2\n_{i}h\n_{\bar{j}}h}{\left(1-h^2\right)^2}|\n h|^2+\frac{2a\n_{i}h\n_{\bar{j}}h}{1-h^2}\left(1+\frac{\ln h}{1-h^2}\right).\label{constrained}
\end{split}
\end{equation}
\end{lemma}
Proof: The proof is straightforward but quite involved, so we divide our proof into fours steps. First, we need to give the evolution equation of $h$. Recalling the definition of $h$, we have
\begin{equation}
\frac{\partial}{\partial t}h=\frac{\partial}{\partial t}\left(\frac{v}{u}\right)=\frac{v_t}{u}-\frac{vu_t}{u^2}=\frac{1}{u^2}\left(u\Delta v-v\Delta u\right)+a\frac{v}{u}\ln\frac{v}{u}\label{i1}
\end{equation}
and
\begin{equation}
\begin{split}
\Delta h=\nabla_k\nabla_{\bar{k}}h=&\nabla_k\frac{u\nabla_{\bar{k}}v-v\nabla_{\bar{k}}u}{u^2}
\\=&\frac{\nabla_ku\nabla_{\bar{k}}v+u\nabla_k\nabla_{\bar{k}}v-\nabla_kv\nabla_{\bar{k}}u-v\nabla_k\nabla_{\bar{k}}u}{u^2}
\\&-\frac{2\nabla_ku\left(u\nabla_{\bar{k}}v-v\nabla_{\bar{k}}u\right)}{u^3}
\\=&\frac{u\Delta v-v\Delta u}{u^2}-\frac{\nabla_{k}u\nabla_{\bar{k}}v+\nabla_kv\nabla_{\bar{k}}u}{u^2}+\frac{2v\nabla_{\bar{k}}u\nabla_{k}u}{u^3}
\\=&\frac{u\Delta v-v\Delta u}{u^2}-\nabla_{k}L\nabla_{\bar{k}}h-\nabla_{\bar{k}}L\nabla_{k}h.\label{i2}
\end{split}
\end{equation}
Combining (\ref{i1}) and (\ref{i2}), we obtain the evolution of $h$
\begin{equation}
\begin{split}
\frac{\partial}{\partial t}h=&\Delta h+\nabla_{k}L\nabla_{\bar{k}}h+\nabla_{\bar{k}}L\nabla_{k}h+ah\ln h.
\end{split}
\end{equation}

Next we will give the evolution equation of $1-h^2$. By
\begin{equation}
\begin{split}
\frac{\partial}{\partial t}\left(1-h^2\right)=&-2h\frac{\partial}{\partial t}h\\=&-2h\left(\Delta h+\nabla_{k}L\nabla_{\bar{k}}h+\nabla_{\bar{k}}L\nabla_{k}h+ah\ln h\right)
\end{split}
\end{equation}
and
\begin{equation}
\Delta\left(1-h^2\right)=-2|\nabla h|^2-2h\Delta h,
\end{equation}
we have that
\begin{equation}
\begin{split}
\frac{\partial}{\partial t}\left(1-h^2\right)=&\Delta\left(1-h^2\right)+2|\nabla h|^2{}-2h\left(\nabla_{k}L\nabla_{\bar{k}}h+\nabla_{\bar{k}}L\nabla_{k}h\right)-2ah^2\ln h\\=&\Delta\left(1-h^2\right)+2|\nabla h|^2+\nabla_{k}L\nabla_{\bar{k}}\left(1-h^2\right)+\nabla_{\bar{k}}L\nabla_{k}\left(1-h^2\right)-2ah^2\ln h. \label{gradient evol1}
\end{split}
\end{equation}

Another step in the proof is computing the evolution equation of $\nabla_{i}h\nabla_{\bar{j}}h$. Take partial derivative with respective to $t$
\begin{equation}
\frac{\partial}{\partial t}\left(\nabla_{i}h\nabla_{\bar{j}}h\right)=\nabla_{i}\left(\frac{\partial}{\partial t}h\right)\nabla_{\bar{j}}h+\nabla_{i}h\nabla_{\bar{j}}\left(\frac{\partial}{\partial t}h\right).
\end{equation}
Direct computation yields
\begin{equation}
\begin{split}
\nabla_{i}\left(\frac{\partial}{\partial t}h\right)\nabla_{\bar{j}}h=&\nabla_{i}\left(\Delta h+\nabla_{k}L\nabla_{\bar{k}}h+\nabla_{\bar{k}}L\nabla_{k}h+ah\ln h\right)\nabla_{\bar{j}}h
\\=&\nabla_{\bar{j}}h\Delta\left(\nabla_{i}h\right)-\frac{1}{2}R_{i\bar{l}}\nabla_{l}h\nabla_{\bar{j}}h+\nabla_i\nabla_{k}L\nabla_{\bar{k}}h\nabla_{\bar{j}}h
\\&+\nabla_{k}L\nabla_i\nabla_{\bar{k}}h\nabla_{\bar{j}}h+\nabla_i\nabla_{\bar{k}}L\nabla_{k}h\nabla_{\bar{j}}h+\nabla_{\bar{k}}L\nabla_i\nabla_{k}h\nabla_{\bar{j}}h\\
&+a(1+\ln h)\nabla_{i}h\nabla_{\bar{j}}h
\end{split}\label{i3}
\end{equation}
and
\begin{equation}
\begin{split}
\nabla_{\bar{j}}\left(\frac{\partial}{\partial t}h\right)\nabla_{i}h=&\nabla_{\bar{j}}\left(\Delta h+\nabla_{k}L\nabla_{\bar{k}}h+\nabla_{\bar{k}}L\nabla_{k}h+ah\ln h\right)\nabla_{i}h
\\=&\n_{\bar{j}}\Delta h\n_{i}h+\n_{\bar{j}}\n_{k}L\n_{\bar{k}}h\n_{i}h+\n_{k}L\n_{\bar{j}}\n_{\bar{k}}h\n_{i}h\\
&+\n_{\bar{j}}\n_{\bar{k}}L\n_{k}h\n_{i}h+\n_{\bar{k}}L\n_{\bar{j}}\n_{k}h\n_{i}h\\
=&\Delta\left(\n_{\bar{j}}h\right)\n_{i}h-\frac{1}{2}R_{k\bar{j}}\n_{\bar{k}}h\n_{i}h\\
&+\n_{\bar{j}}\n_{k}L\n_{\bar{k}}h\n_{i}h+\n_{k}L\n_{\bar{j}}\n_{\bar{k}}h\n_{i}h\\
&+\n_{\bar{j}}\n_{\bar{k}}L\n_{k}h\n_{i}h+\n_{\bar{k}}L\n_{\bar{j}}\n_{k}h\n_{i}h\\
&+a(1+\ln h)\nabla_{i}h\nabla_{\bar{j}}h.
\end{split}\label{i4}
\end{equation}
Combining (\ref{i3}) and (\ref{i4}), we have that
\begin{equation}
\begin{split}
\frac{\partial}{\partial t}\left(\n_{i}h\n_{\bar{j}}h\right)=&\Delta\left(\n_{i}h\right)\n_{\bar{j}}h+\Delta\left(\n_{\bar{j}}h\right)\n_{i}h-\frac{1}{2}R_{i\bar{l}}\n_{l}h\n_{\bar{j}}h\\
&-\frac{1}{2}R_{k\bar{j}}\n_{\bar{k}}h\n_{i}h+\n_{i}\n_{k}L\n_{\bar{j}}h\n_{\bar{k}}h+\n_{i}\n_{\bar{k}}L\n_{k}h\n_{\bar{j}}h\\
&+\n_{k}L\n_{\bar{k}}\left(\n_{i}h\n_{\bar{j}}h\right)+\n_{\bar{k}}L\n_{k}\left(\n_{i}h\n_{\bar{j}}h\right)
\\&+\n_{\bar{j}}\n_{k}L\n_{\bar{k}}h\n_{i}h+\n_{\bar{j}}\n_{\bar{k}}L\n_{k}h\n_{i}h+2a(1+\ln h)\nabla_{i}h\nabla_{\bar{j}}h.
\end{split}
\end{equation}
The Laplacian of $\nabla_{i}h\nabla_{\bar{j}}h$ is
\begin{equation}
\begin{split}
\Delta\left(\n_{i}h\n_{\bar{j}}h\right)
=&\frac{1}{2}\left(\n_{k}\n_{\bar{k}}+\n_{\bar{k}}\n_{k}\right)\left(\n_{i}h\n_{\bar{j}}h\right)\\
=&\frac{1}{2}\n_{k}\left(\n_{\bar{k}}\n_{i}h\n_{\bar{j}}h+\n_{i}h\n_{\bar{k}}\n_{\bar{j}}h\right)
\\&+\frac{1}{2}\n_{\bar{k}}\left(\n_{k}\n_{i}h\n_{\bar{j}}h+\n_{i}h\n_{k}\n_{\bar{j}}h\right)\\
=&\Delta\left(\n_{i}h\right)\n_{\bar{j}}h+\n_{i}h\Delta\left(\n_{\bar{j}}h\right)\\&+\n_{k}\n_{i}h\n_{\bar{k}}\n_{\bar{j}}h+\n_{\bar{k}}\n_{i}h\n_{k}\n_{\bar{j}}h.
\end{split}
\end{equation}
Collecting the two equalities above gives that
\begin{equation}
\begin{split}
\frac{\partial}{\partial t}\left(\n_{i}h\n_{\bar{j}}h\right)=&\Delta\left(\n_{i}h\n_{\bar{j}}h\right)-\n_{k}\n_{i}h\n_{\bar{k}}\n_{\bar{j}}h-\n_{\bar{k}}\n_{i}h\n_{k}\n_{\bar{j}}h
\\&+\n_{i}\n_{k}L\n_{\bar{j}}h\n_{\bar{k}}h
+\n_{i}\n_{\bar{k}}L\n_{k}h\n_{\bar{j}}h\\&+\n_{\bar{j}}\n_{k}L\n_{\bar{k}}h\n_{i}h+\n_{\bar{j}}\n_{\bar{k}}L\n_{k}h\n_{i}h\\
&+\n_{k}L\n_{\bar{k}}\left(\n_{i}h\n_{\bar{j}}h\right)+\n_{\bar{k}}L\n_{k}\left(\n_{i}h\n_{\bar{j}}h\right)\\&
-\frac{1}{2}R_{i\bar{l}}\n_{\bar{j}}h\n_{l}h-\frac{1}{2}R_{k\bar{j}}\n_{\bar{k}}h\n_{i}h+2a(1+\ln h)\nabla_{i}h\nabla_{\bar{j}}h. \label{gradient evol2}
\end{split}
\end{equation}

Finally, we give the evolution equation of $\frac{\n_{i}h\n_{\bar{j}}h}{1-h^2}$. With the help of (\ref{gradient evol1}) and (\ref{gradient evol2}), we have that
\begin{equation}
\begin{split}
\frac{\partial}{\partial t}\left(\frac{\n_{i}h\n_{\bar{j}}h}{1-h^2}\right)=&\frac{1}{1-h^2}\frac{\partial}{\partial t}{\left(\n_{i}h\n_{\bar{j}}h\right)}-\frac{\n_{i}h\n_{\bar{j}}h}{\left(1-h^2\right)^2}\frac{\partial}{\partial t}\left(1-h^2\right)\\=&\frac{1}{1-h^2}[\Delta\left(\n_{i}h\n_{\bar{j}}h\right)-\n_{k}\n_{i}h\n_{\bar{k}}\n_{\bar{j}}h-\n_{\bar{k}}\n_{i}h\n_{k}\n_{\bar{j}}h{}\\&+\n_{i}\n_{k}L\n_{\bar{j}}h\n_{\bar{k}}h
+\n_{i}\n_{\bar{k}}L\n_{k}h\n_{\bar{j}}h+\n_{\bar{j}}\n_{k}L\n_{\bar{k}}h\n_{i}h\\&+\n_{\bar{j}}\n_{\bar{k}}L\n_{i}h\n_{k}h+\n_{k}L\n_{\bar{k}}(\n_{i}h\n_{\bar{j}}h)+\n_{\bar{k}}L\n_{k}\left(\n_{i}h\n_{\bar{j}}h\right)\\
&-\frac{1}{2}R_{i\bar{l}}\n_{\bar{j}}h\n_{l}h-\frac{1}{2}R_{k\bar{j}}\n_{\bar{k}}h\n_{i}h+2a(1+\ln h)\nabla_{i}h\nabla_{\bar{j}}h]\\
&-\frac{\n_{i}h\n_{\bar{j}}h}{\left(1-h^2\right)^2}[\Delta\left(1-h^2\right)+2|\n h|^2+\n_{k}L\n_{\bar{k}}\left(1-h^2\right) \\&+\n_{\bar{k}}L\n_{k}\left(1-h^2\right)+2ah^2\ln h].
\end{split}
\end{equation}
Straightforward computation yields
\begin{equation}
\begin{split}
\Delta\left(\frac{\n_{i}h\n_{\bar{j}}h}{1-h^2}\right)=&\frac{\Delta\left(\n_{i}h\n_{\bar{j}}h\right)}{1-h^2}-\frac{\n_{i}h\n_{\bar{j}}h}{\left(1-h^2\right)^2}\Delta\left(1-h^2\right)-\frac{\n_{k}\left(\n_{i}h\n_{\bar{j}}h\right)\n_{\bar{k}}\left(1-h^2\right)}{\left(1-h^2\right)^2}\\&
-\frac{\n_{\bar{k}}\left(\n_{i}h\n_{\bar{j}}h\right)\n_{k}\left(1-h^2\right)}{\left(1-h^2\right)^2}+\frac{2\n_{i}h\n_{\bar{j}}h[\n_{k}\left(1-h^2\right)\n_{\bar{k}}\left(1-h^2\right)]}{\left(1-h^2\right)^3}\\
=&\frac{\Delta\left(\n_{i}h\n_{\bar{j}}h\right)}{1-h^2}-\frac{\n_{i}h\n_{\bar{j}}h}{\left(1-h^2\right)^2}\Delta\left(1-h^2\right)\\
&+\frac{2h\n_{k}\left(\n_{i}h\n_{\bar{j}}h\right)\n_{\bar{k}}h+2h\n_{\bar{k}}\left(\n_{i}h\n_{\bar{j}}h\right)\n_{k}h}{\left(1-h^2\right)^2}\\&+\frac{8h^2\n_{i}h\n_{\bar{j}}h|\n h|^2}{\left(1-h^2\right)^3},
\end{split}
\end{equation}
so we have that
\begin{equation}
\begin{split}
\frac{\partial}{\partial t}\left(\frac{\n_{i}h\n_{\bar{j}}h}{1-h^2}\right)
=&\Delta\left(\frac{\n_{i}h\n_{\bar{j}}h}{1-h^2}\right)+\frac{1}{1-h^2}(-\n_{k}\n_{i}h\n_{\bar{k}}\n_{\bar{j}}h\\&-\n_{\bar{k}}\n_{i}h\n_{k}\n_{\bar{j}}h+\n_{i}\n_{k}L\n_{\bar{j}}h\n_{\bar{k}}h\\
&+\n_{i}\n_{\bar{k}}L\n_{\bar{j}}h\n_{k}h+\n_{\bar{j}}\n_{k}L\n_{\bar{k}}h\n_{i}h\\
&+\n_{\bar{j}}\n_{\bar{k}}L\n_{i}h\n_{k}h-\frac{1}{2}R_{i\bar{l}}\n_{\bar{j}}h\n_{l}h-\frac{1}{2}R_{k\bar{j}}\n_{\bar{k}}h\n_{i}h)
\\&-\frac{2\n_{i}h\n_{\bar{j}}h}{\left(1-h^2\right)^2}|\n h|^2+\n_{k}L\n_{\bar{k}}\left(\frac{\n_{i}h\n_{\bar{j}}h}{1-h^2}\right)\\
&+\n_{\bar{k}}L\n_{k}\left(\frac{\n_{i}h\n_{\bar{j}}h}{1-h^2}\right)-\frac{8h^2\n_{i}h\n_{\bar{j}}h|\n h|^2}{\left(1-h^2\right)^3}\\
&-\frac{2h\n_{k}\left(\n_{i}h\n_{\bar{j}}h\right)\n_{\bar{k}}h+2h\n_{\bar{k}}\left(\n_{i}h\n_{\bar{j}}h\right)\n_{k}h}{\left(1-h^2\right)^2}\\
&+\frac{2a\n_{i}h\n_{\bar{j}}h}{1-h^2}\left(1+\ln h+\frac{h^2\ln h}{1-h^2}\right).
\end{split}
\end{equation}
Rearranging terms leads to
\begin{equation}
\begin{split}
\frac{\partial}{\partial t}\left(\frac{\n_{i}h\n_{\bar{j}}h}{1-h^2}\right)
=&\Delta\left(\frac{\n_{i}h\n_{\bar{j}}h}{1-h^2}\right)+\n_{k}L\n_{\bar{k}}\left(\frac{\n_{i}h\n_{\bar{j}}h}{1-h^2}\right)+\n_{\bar{k}}L\n_{k}\left(\frac{\n_{i}h\n_{\bar{j}}h}{1-h^2}\right)\\
&-\frac{1}{1-h^2}\left(\n_{i}\n_{k}h+\frac{2h\n_{i}h\n_{k}h}{1-h^2}\right)\left(\n_{\bar{j}}\n_{\bar{k}}h+\frac{2h\n_{\bar{j}}h\n_{\bar{k}}h}{1-h^2}\right)\\
&-\frac{1}{1-h^2}\left(\n_{i}\n_{\bar{k}}h+\frac{2h\n_{i}h\n_{\bar{k}}h}{1-h^2}\right)\left(\n_{\bar{j}}\n_{k}h+\frac{2h\n_{\bar{j}}h\n_{k}h}{1-h^2}\right)\\
&+\frac{1}{1-h^2}(\n_{i}\n_{k}L\n_{\bar{j}}h\n_{\bar{k}}h+\n_{i}\n_{\bar{k}}L\n_{\bar{j}}h\n_{k}h+\n_{\bar{j}}\n_{k}L\n_{i}h\n_{\bar{k}}h\\&+\n_{\bar{j}}\n_{\bar{k}}L\n_{i}h\n_{k}h
-\frac{1}{2}R_{i\bar{l}}\n_{\bar{j}}h\n_{l}h-\frac{1}{2}R_{k\bar{j}}\n_{\bar{k}}h\n_{i}h)\\
&-\frac{2\n_{i}h\n_{\bar{j}}h}{\left(1-h^2\right)^2}|\n h|^2+\frac{2a\n_{i}h\n_{\bar{j}}h}{1-h^2}\left(1+\frac{\ln h}{1-h^2}\right).
\end{split}
\end{equation}
This completes the proof of Lemma 4.1.

Similar calculation leads to
\begin{lemma}If $u$ and $v$ are positive solutions to the nonlinear heat equation (\ref{NLF}), then the radio $h=v/u$ also satisfies (\ref{constrained}).
\end{lemma}

Now we are ready to give the proof of the Theorem 1.3 and Theorem 1.4.

{\bf Proof of the Theorem 1.3}
Let \begin{equation}
P_{i{\bar{j}}}=\n_{i}\n_{\bar{j}}L-\frac{\n_{i}h\n_{\bar{j}}h}{1-h^2}.
\end{equation}
Combine the lemmas together to get
\begin{equation}
\begin{split}
\frac{\partial}{\partial t}P_{i{\bar{j}}}=&\frac{\partial}{\partial t}\n_{i}\n_{\bar{j}}L-\frac{\partial}{\partial t}\left(\frac{\n_{i}h\n_{\bar{j}}h}{1-h^2}\right)\\
=&\Delta\left(\nabla_i\nabla_{\bar{j}}L\right)+R_{i\bar{j}l\bar{k}}\nabla_k\nabla_{\bar{l}}L
+R_{i\bar{j}k\bar{l}}{\nabla_{l}L\nabla_{\bar{k}}L}
+\nabla_{k}L\nabla_{\bar{k}}\nabla_i\nabla_{\bar{j}}L\\
&+\nabla_{\bar{k}}L\nabla_k\nabla_i\nabla_{\bar{j}}L
+\nabla_{\bar{j}}\nabla_{k}L\nabla_i\nabla_{\bar{k}}L+\nabla_i\nabla_{k}L\nabla_{\bar{j}}\nabla_{\bar{k}}L
\\&-\frac{1}{2}\left(R_{l\bar{j}}\nabla_i\nabla_{\bar{l}}L+R_{i\bar{l}}\nabla_l\nabla_{\bar{j}}L\right)+a\nabla_i\nabla_{\bar{j}}L\\
&-\Delta\left(\frac{\n_{i}h\n_{\bar{j}}h}{1-h^2}\right)-\n_{k}L\n_{\bar{k}}\left(\frac{\n_{i}h\n_{\bar{j}}h}{1-h^2}\right)-\n_{\bar{k}}L\n_{k}\left(\frac{\n_{i}h\n_{\bar{j}}h}{1-h^2}\right)\\
&+\frac{1}{1-h^2}\left(\n_{i}\n_{k}h+\frac{2h\n_{i}h\n_{k}h}{1-h^2}\right)\left(\n_{\bar{j}}\n_{\bar{k}}h+\frac{2h\n_{\bar{j}}h\n_{\bar{k}}h}{1-h^2}\right)\\
&+\frac{1}{1-h^2}\left(\n_{i}\n_{\bar{k}}h+\frac{2h\n_{i}h\n_{\bar{k}}h}{1-h^2}\right)\left(\n_{\bar{j}}\n_{k}h+\frac{2h\n_{\bar{j}}h\n_{k}h}{1-h^2}\right)\\
&-\n_{i}\n_{k}L\left(\frac{\n_{\bar{j}}h\n_{\bar{k}}h}{1-h^2}\right)-\n_{i}\n_{\bar{k}}L\left(\frac{\n_{\bar{j}}h\n_{k}h}{1-h^2}\right)-\n_{\bar{j}}\n_{k}L\left(\frac{\n_{i}h\n_{\bar{k}}h}{1-h^2}\right)
\\&-\n_{\bar{j}}\n_{\bar{k}}L\left(\frac{\n_{i}h\n_{k}h}{1-h^2}\right)+\frac{1}{2}R_{i\bar{l}}\left(\frac{\n_{\bar{j}}h\n_{l}h}{1-h^2}\right)+\frac{1}{2}R_{k\bar{j}}\left(\frac{\n_{\bar{k}}h\n_{i}h}{1-h^2}\right)\\
&+\frac{2\n_{i}h\n_{\bar{j}}h}{\left(1-h^2\right)^2}|\n h|^2-\frac{2a\n_{i}h\n_{\bar{j}}h}{1-h^2}\left(1+\frac{\ln h}{1-h^2}\right)\\
=&\Delta P_{i{\bar{j}}} +\nabla_{k}L\nabla_{\bar{k}}P_{i{\bar{j}}}+\nabla_{\bar{k}}L\nabla_{k}P_{i{\bar{j}}}+R_{i{\bar{j}}k{\bar{l}}}P_{l{\bar{k}}}+R_{i\bar{j}k\bar{l}}{\nabla_{l}L\nabla_{\bar{k}}L}+aP_{i\bar{j}}
\\&-\frac{1}{2}\left(R_{l{\bar{j}}}P_{i{\bar{l}}}+R_{i{\bar{l}}}P_{l{\bar{j}}}\right)+R_{i{\bar{j}}k{\bar{l}}}\frac{\n_{l}h\n_{\bar{k}}h}{1-h^2}-\frac{a\n_{i}h\n_{\bar{j}}h}{1-h^2}\left(1+\frac{2\ln h}{1-h^2}\right)\\
&+\left(\nabla_{\bar{j}}\nabla_{k}L-\frac{\n_{\bar{j}}h\n_{k}h}{1-h^2}\right)\left(\nabla_{i}\nabla_{\bar{k}}L-\frac{\n_{i}h\n_{\bar{k}}h}{1-h^2}\right)\\
&+\left(\n_{i}\n_{k}L-\frac{\n_{i}h\n_{k}h}{1-h^2}\right)\left(\n_{\bar{j}}\n_{\bar{k}}L-\frac{\n_{\bar{j}}h\n_{\bar{k}}h}{1-h^2}\right)\\&+\frac{1}{1-h^2}\left(\n_{i}\n_{k}h+\frac{2h\n_{i}h\n_{k}h}{1-h^2}\right)\left(\n_{\bar{j}}\n_{\bar{k}}h+\frac{2h\n_{\bar{j}}h\n_{\bar{k}}h}{1-h^2}\right)\\
&+\frac{1}{1-h^2}\left(\n_{i}\n_{\bar{k}}h+\frac{2h\n_{i}h\n_{\bar{k}}h}{1-h^2}\right)\left(\n_{\bar{j}}\n_{k}h+\frac{2h\n_{\bar{j}}h\n_{k}h}{1-h^2}\right).
\end{split}
\end{equation}
Since the three matrices in the last three terms are nonnegative definite, we have that
\begin{equation}
\begin{split}
\frac{\partial}{\partial t}P_{i{\bar{j}}}\geq&\Delta P_{i{\bar{j}}} +\nabla_{k}L\nabla_{\bar{k}}P_{i{\bar{j}}}+\nabla_{\bar{k}}L\nabla_{k}P_{i{\bar{j}}}+R_{i{\bar{j}}k{\bar{l}}}P_{l{\bar{k}}}+R_{i\bar{j}k\bar{l}}{\nabla_{l}L\nabla_{\bar{k}}L}+aP_{i\bar{j}}
\\&+P_{i\bar{k}}P_{k\bar{j}}-\frac{1}{2}\left(R_{l{\bar{j}}}P_{i{\bar{l}}}+R_{i{\bar{l}}}P_{l{\bar{j}}}\right)
+R_{i{\bar{j}}k{\bar{l}}}\frac{\n_{l}h\n_{\bar{k}}h}{1-h^2}-\frac{a\n_{i}h\n_{\bar{j}}h}{1-h^2}\left(1+\frac{2\ln h}{1-h^2}\right)
\end{split}
\end{equation}
It was shown in \cite{W} that the function
$$G(h)=1+\frac{2\ln h}{1-h^2}$$
is increasing and negative for $h\in (0,1)$.\\
If $a>0$, then by the nonnegativity of holomorphic bisectional curvature we have
$$
\frac{\partial}{\partial t}P_{i{\bar{j}}}\geq\Delta P_{i{\bar{j}}} +\nabla_{k}L\nabla_{\bar{k}}P_{i{\bar{j}}}+\nabla_{\bar{k}}L\nabla_{k}P_{i{\bar{j}}}+R_{i{\bar{j}}k{\bar{l}}}P_{l{\bar{k}}}+aP_{i\bar{j}}
+P_{i\bar{k}}P_{k\bar{j}}-\frac{1}{2}\left(R_{l{\bar{j}}}P_{i{\bar{l}}}+R_{i{\bar{l}}}P_{l{\bar{j}}}\right)
$$
Now letting
$$Q_{i\bar{j}}=P_{i\bar{j}}+\frac{a}{1-e^{-at}}g_{i\bar{j}},$$
we have
\begin{equation}
\begin{split}
\frac{\partial}{\partial t}Q_{i\bar{j}}=&\frac{\partial}{\partial t}\left(P_{i{\bar{j}}}+\frac{a}{1-e^{-at}}g_{i\bar{j}}\right)\\
=&\frac{\partial}{\partial t}P_{i{\bar{j}}}-\frac{a^2e^{-at}}{(1-e^{-at})^2}g_{i{\bar{j}}}\\
\geq &\Delta P_{i{\bar{j}}} +\nabla_{k}L\nabla_{\bar{k}}P_{i{\bar{j}}}+\nabla_{\bar{k}}L\nabla_{k}P_{i{\bar{j}}}+R_{i{\bar{j}}k{\bar{l}}}P_{l{\bar{k}}}+aP_{i\bar{j}}
\\&+P_{i\bar{k}}P_{k\bar{j}}-\frac{1}{2}\left(R_{l{\bar{j}}}P_{i{\bar{l}}}+R_{i{\bar{l}}}P_{l{\bar{j}}}\right)-\frac{a^2e^{-at}}{(1-e^{-at})^2}g_{i{\bar{j}}}\\
=&\Delta Q_{i{\bar{j}}} +\nabla_{k}L\nabla_{\bar{k}}Q_{i{\bar{j}}}+\nabla_{\bar{k}}L\nabla_{k}Q_{i{\bar{j}}}+R_{i{\bar{j}}k{\bar{l}}}Q_{l\bar{k}}+aQ_{i\bar{j}}-\frac{a^2}{1-e^{-at}}g_{i\bar{j}}
\\&+(Q_{i\bar{k}}-\frac{a}{1-e^{-at}}g_{i\bar{k}})(Q_{k\bar{j}}-\frac{a}{1-e^{-at}}g_{k\bar{j}})-\frac{1}{2}\left(R_{l{\bar{j}}}Q_{i{\bar{l}}}+R_{i{\bar{l}}}Q_{l{\bar{j}}}\right)\\
&-\frac{a^2e^{-at}}{(1-e^{-at})^2}g_{i{\bar{j}}}\\
=&\Delta Q_{i{\bar{j}}} +\nabla_{k}L\nabla_{\bar{k}}Q_{i{\bar{j}}}+\nabla_{\bar{k}}L\nabla_{k}Q_{i{\bar{j}}}+R_{i{\bar{j}}k{\bar{l}}}Q_{l\bar{k}}
+Q_{i\bar{k}}Q_{k\bar{j}}-\frac{a(1+e^{-at})}{1-e^{-at}}Q_{i\bar{j}}\\
&-\frac{1}{2}\left(R_{l{\bar{j}}}Q_{i{\bar{l}}}+R_{i{\bar{l}}}Q_{l{\bar{j}}}\right)\end{split}
\end{equation}
It follows from Hamilton's maximum principle for tensor that $Q_{i{\bar{j}}}\geq 0$.
\newline Next we discuss the case $a<0$. It follows from the assumption of curvature tensor that
\begin{equation}
R_{i{\bar{j}}k{\bar{l}}}\frac{\n_{l}h\n_{\bar{k}}h}{1-h^2}\geq \frac{-aK}{1-h^2}(g_{i{\bar{j}}}|\nabla h|^2+\n_{i}h\n_{\bar{j}}h)\geq\frac{-aK|\nabla h|^2}{1-h^2}g_{i{\bar{j}}}.\label{4.22}
\end{equation}
Thanks to the negativity of $G(h)$, we have that
\begin{equation}
-\frac{a\n_{i}h\n_{\bar{j}}h}{1-h^2}\left(1+\frac{2\ln h}{1-h^2}\right)\geq-\frac{a|\nabla h|^2}{1-h^2}\left(1+\frac{2\ln h}{1-h^2}\right)g_{i{\bar{j}}}.\label{4.23}
\end{equation}
Combining (\ref{4.22}) and (\ref{4.23}), we obtain
\begin{equation}
R_{i{\bar{j}}k{\bar{l}}}\frac{\n_{l}h\n_{\bar{k}}h}{1-h^2}-\frac{a\n_{i}h\n_{\bar{j}}h}{1-h^2}\left(1+\frac{2\ln h}{1-h^2}\right)\geq -\frac{a|\nabla h|^2}{1-h^2}\left(K+1+\frac{2\ln h}{1-h^2}\right)g_{i{\bar{j}}}.
\end{equation}
Since $G(h)$ is increasing and $K$ is defined such that
$$K+1+\frac{2\ln h}{1-h^2}\geq 0,$$
 we get
$$
\frac{\partial}{\partial t}P_{i{\bar{j}}}\geq\Delta P_{i{\bar{j}}} +\nabla_{k}L\nabla_{\bar{k}}P_{i{\bar{j}}}+\nabla_{\bar{k}}L\nabla_{k}P_{i{\bar{j}}}+R_{i{\bar{j}}k{\bar{l}}}P_{l{\bar{k}}}+aP_{i\bar{j}}
+P_{i\bar{k}}P_{k\bar{j}}-\frac{1}{2}\left(R_{l{\bar{j}}}P_{i{\bar{l}}}+R_{i{\bar{l}}}P_{l{\bar{j}}}\right).
$$
Now we can easily carry out the proof of this case by using the same argument as in the proof of the case above.

{\bf Proof of the Theorem 1.4}
Let \begin{equation}
P_{i{\bar{j}}}=\n_{i}\n_{\bar{j}}L+R_{i\bar{j}}-\frac{\n_{i}h\n_{\bar{j}}h}{1-h^2}.
\end{equation}
Combine Lemma 3.5 and Lemma 4.2 together to get
\begin{equation}
\begin{split}
\frac{\partial}{\partial t}P_{i{\bar{j}}}=&\frac{\partial}{\partial t}\n_{i}\n_{\bar{j}}L+\frac{\partial}{\partial t}R_{i\bar{j}}-\frac{\partial}{\partial t}\left(\frac{\n_{i}h\n_{\bar{j}}h}{1-h^2}\right)\\
=&\Delta\left(\nabla_i\nabla_{\bar{j}}L\right)+R_{i\bar{j}l\bar{k}}\nabla_k\nabla_{\bar{l}}L
+R_{i\bar{j}k\bar{l}}{\nabla_{l}L\nabla_{\bar{k}}L}
+\nabla_{k}L\nabla_{\bar{k}}\nabla_i\nabla_{\bar{j}}L\\
&+\nabla_{\bar{k}}L\nabla_k\nabla_i\nabla_{\bar{j}}L
+\nabla_{\bar{j}}\nabla_{k}L\nabla_i\nabla_{\bar{k}}L+\nabla_i\nabla_{k}L\nabla_{\bar{j}}\nabla_{\bar{k}}L
\\&-\frac{1}{2}\left(R_{l\bar{j}}\nabla_i\nabla_{\bar{l}}L+R_{i\bar{l}}\nabla_l\nabla_{\bar{j}}L\right)+a\nabla_i\nabla_{\bar{j}}L\\
&+2\Delta R_{i{\bar{j}}}+2R_{i{\bar{j}}k{\bar{l}}}R_{l{\bar{k}}}-2R_{i{\bar{p}}}R_{p{\bar{j}}}\\
&-\Delta\left(\frac{\n_{i}h\n_{\bar{j}}h}{1-h^2}\right)-\n_{k}L\n_{\bar{k}}\left(\frac{\n_{i}h\n_{\bar{j}}h}{1-h^2}\right)-\n_{\bar{k}}L\n_{k}\left(\frac{\n_{i}h\n_{\bar{j}}h}{1-h^2}\right)\\
&+\frac{1}{1-h^2}\left(\n_{i}\n_{k}h+\frac{2h\n_{i}h\n_{k}h}{1-h^2}\right)\left(\n_{\bar{j}}\n_{\bar{k}}h+\frac{2h\n_{\bar{j}}h\n_{\bar{k}}h}{1-h^2}\right)\\
&+\frac{1}{1-h^2}\left(\n_{i}\n_{\bar{k}}h+\frac{2h\n_{i}h\n_{\bar{k}}h}{1-h^2}\right)\left(\n_{\bar{j}}\n_{k}h+\frac{2h\n_{\bar{j}}h\n_{k}h}{1-h^2}\right)\\
&-\n_{i}\n_{k}L\left(\frac{\n_{\bar{j}}h\n_{\bar{k}}h}{1-h^2}\right)-\n_{i}\n_{\bar{k}}L\left(\frac{\n_{\bar{j}}h\n_{k}h}{1-h^2}\right)-\n_{\bar{j}}\n_{k}L\left(\frac{\n_{i}h\n_{\bar{k}}h}{1-h^2}\right)
\\&-\n_{\bar{j}}\n_{\bar{k}}L\left(\frac{\n_{i}h\n_{k}h}{1-h^2}\right)+\frac{1}{2}R_{i\bar{l}}\left(\frac{\n_{\bar{j}}h\n_{l}h}{1-h^2}\right)+\frac{1}{2}R_{k\bar{j}}\left(\frac{\n_{\bar{k}}h\n_{i}h}{1-h^2}\right)\\
&+\frac{2\n_{i}h\n_{\bar{j}}h}{\left(1-h^2\right)^2}|\n h|^2-\frac{2a\n_{i}h\n_{\bar{j}}h}{1-h^2}\left(1+\frac{\ln h}{1-h^2}\right)\\
=&\Delta P_{i{\bar{j}}} +\nabla_{k}L\nabla_{\bar{k}}P_{i{\bar{j}}}+\nabla_{\bar{k}}L\nabla_{k}P_{i{\bar{j}}}+R_{i{\bar{j}}k{\bar{l}}}P_{l{\bar{k}}}+aP_{i\bar{j}}-aR_{i{\bar{j}}}\\
&-\frac{1}{2}\left(R_{l{\bar{j}}}P_{i{\bar{l}}}+R_{i{\bar{l}}}P_{l{\bar{j}}}\right)-R_{i{\bar{p}}}R_{p{\bar{j}}}+R_{i{\bar{j}}k{\bar{l}}}\frac{\n_{l}h\n_{\bar{k}}h}{1-h^2}-\frac{a\n_{i}h\n_{\bar{j}}h}{1-h^2}\left(1+\frac{2\ln h}{1-h^2}\right)\\
&+\Delta R_{i{\bar{j}}}+R_{i{\bar{j}}k{\bar{l}}}R_{l{\bar{k}}}-\nabla_{k}L\nabla_{\bar{k}}R_{i{\bar{j}}}-\nabla_{\bar{k}}L\nabla_kR_{i{\bar{j}}}+R_{i\bar{j}k\bar{l}}{\nabla_{l}L\nabla_{\bar{k}}L}\\
&+\left(\nabla_{\bar{j}}\nabla_{k}L-\frac{\n_{\bar{j}}h\n_{k}h}{1-h^2}\right)\left(\nabla_{i}\nabla_{\bar{k}}L-\frac{\n_{i}h\n_{\bar{k}}h}{1-h^2}\right)\\
&+\left(\n_{i}\n_{k}L-\frac{\n_{i}h\n_{k}h}{1-h^2}\right)\left(\n_{\bar{j}}\n_{\bar{k}}L-\frac{\n_{\bar{j}}h\n_{\bar{k}}h}{1-h^2}\right)\\&+\frac{1}{1-h^2}\left(\n_{i}\n_{k}h+\frac{2h\n_{i}h\n_{k}h}{1-h^2}\right)\left(\n_{\bar{j}}\n_{\bar{k}}h+\frac{2h\n_{\bar{j}}h\n_{\bar{k}}h}{1-h^2}\right)\\
&+\frac{1}{1-h^2}\left(\n_{i}\n_{\bar{k}}h+\frac{2h\n_{i}h\n_{\bar{k}}h}{1-h^2}\right)\left(\n_{\bar{j}}\n_{k}h+\frac{2h\n_{\bar{j}}h\n_{k}h}{1-h^2}\right).
\end{split}
\end{equation}
Since the three matrices in the last three terms are nonnegative definite, we have that
\begin{equation}
\begin{split}
\frac{\partial}{\partial t}P_{i{\bar{j}}}\geq&\Delta P_{i{\bar{j}}} +\nabla_{k}L\nabla_{\bar{k}}P_{i{\bar{j}}}+\nabla_{\bar{k}}L\nabla_{k}P_{i{\bar{j}}}+R_{i{\bar{j}}k{\bar{l}}}P_{l{\bar{k}}}+aP_{i\bar{j}}+(P_{i{\bar k}}-R_{i{\bar k}})(P_{k{\bar j}}-R_{k{\bar j}})\\
&-\frac{1}{2}\left(R_{l{\bar{j}}}P_{i{\bar{l}}}+R_{i{\bar{l}}}P_{l{\bar{j}}}\right)-R_{i{\bar{p}}}R_{p{\bar{j}}}+R_{i{\bar{j}}k{\bar{l}}}\frac{\n_{l}h\n_{\bar{k}}h}{1-h^2}-\frac{a\n_{i}h\n_{\bar{j}}h}{1-h^2}\left(1+\frac{2\ln h}{1-h^2}\right)\\
&+\Delta R_{i{\bar{j}}}+R_{i{\bar{j}}k{\bar{l}}}R_{l{\bar{k}}}-\nabla_{k}L\nabla_{\bar{k}}R_{i{\bar{j}}}-\nabla_{\bar{k}}L\nabla_kR_{i{\bar{j}}}+R_{i\bar{j}k\bar{l}}{\nabla_{l}L\nabla_{\bar{k}}L}-aR_{i{\bar{j}}}.
\end{split}
\end{equation}
By Lemma 3.4 and the nonnegativity of holomorphic bisectional curvature, we have
\begin{equation}
\begin{split}
\frac{\partial}{\partial t}P_{i{\bar{j}}}\geq&\Delta P_{i{\bar{j}}} +\nabla_{k}L\nabla_{\bar{k}}P_{i{\bar{j}}}+\nabla_{\bar{k}}L\nabla_{k}P_{i{\bar{j}}}+R_{i{\bar{j}}k{\bar{l}}}P_{l{\bar{k}}}+P_{i{\bar k}}P_{k{\bar j}}\\
&-\frac{3}{2}\left(R_{l{\bar{j}}}P_{i{\bar{l}}}+R_{i{\bar{l}}}P_{l{\bar{j}}}\right)+aP_{i\bar{j}}-aR_{i\bar{j}}-\frac{ae^{-at}}{1-e^{-at}}R_{i\bar{j}}
\end{split}
\end{equation}
Now the theorem can be proved in the same way as shown in the proof of Theorem 1.2.

\section{Extensions}
\setcounter{equation}{0}

In this section, we consider the matrix Li-Yau-Hamilton estimates for more general nonlinear heat equation. The following lemma, which was proven in \cite{CH}, is important in the proof of Theorem 1.5.
\begin{lemma}
If a function $A$ satisfies the equation
\begin{equation}
\frac{\partial}{\partial t}A=\Delta A+B,
\end{equation}
then
\begin{equation}
\begin{split}
\frac{\partial}{\partial t}\nabla_i\nabla_i A
=&\Delta(\nabla_i\nabla_i A)+2R_{kijl}\nabla_k\nabla_l B-R_{il}\nabla_j\nabla_lB-R_{jl}\nabla_i\nabla_lB\\
&-(\nabla_iR_{jl}+\nabla_jR_{il}-\nabla_lR_{ij})\nabla_lB+\nabla_i\nabla_lB.
\end{split}
\end{equation}
\end{lemma}

Now  we are going to prove Theorem 1.5 and Theorem 1.6.

{\bf Proof of Theorem 1.5} Set $$P_{ij}=\nabla_i\nabla_j L+f(t)g_{ij}$$ and we have that
\begin{equation}
\begin{split}
\frac{\partial}{\partial t}P_{ij}=&\frac{\partial}{\partial t}\nabla_i\nabla_j L+f'(t)g_{ij}\\
=&\Delta(\nabla_i\nabla_i L)+2R_{kijl}\nabla_k\nabla_l L-R_{il}\nabla_j\nabla_lL-R_{jl}\nabla_i\nabla_lL\\
&-(\nabla_iR_{jl}+\nabla_jR_{il}-\nabla_lR_{ij})\nabla_lL+2\nabla_i\nabla_lL\cdot\nabla_j\nabla_lL+2\nabla_lL\cdot\nabla_l\nabla_i\nabla_j L\\
&+2R_{kijl}\nabla_k L\nabla_l L+F'(L)\nabla_i\nabla_j L+F''(L)\nabla_iL\nabla_j L+f'(t)g_{ij}\\
=&\Delta P_{ij}+2R_{kijl}P_{kl}-R_{il}P_{jl}-R_{jl}P_{il}-(\nabla_iR_{jl}+\nabla_jR_{il}-\nabla_lR_{ij})\nabla_lL\\
&+2P_{il}P_{jl}-4f(t)P_{ij}+2f^2(t)g_{ij}+2\nabla_lL\cdot\nabla_lP_{ij}+2R_{kijl}\nabla_k L\nabla_l L\\
&+F'(L)P_{ij}-F'(L)f(t)g_{ij}+F''(L)\nabla_iL\nabla_j L+f'(t)g_{ij}\\
=&\Delta P_{ij}+2R_{kijl}P_{kl}-R_{il}P_{jl}-R_{jl}P_{il}-(\nabla_iR_{jl}+\nabla_jR_{il}-\nabla_lR_{ij})\nabla_lL\\
&+2P_{il}P_{jl}+(F'(L)-4f(t))P_{ij}+2\nabla_lL\cdot\nabla_lP_{ij}+2R_{ikjl}\nabla_k L\nabla_l L\\
&+F''(L)\nabla_i L\nabla_j L+(2f^2(t)-F'(L)f(t)+f'(t))g_{ij}\\
\geq&\Delta P_{ij}+2\nabla_lL\cdot\nabla_lP_{ij}+2R_{kijl}P_{kl}-R_{il}P_{jl}-R_{jl}P_{il}+2P_{il}P_{jl}+F''(L)\nabla_i L\nabla_j L\\
&+(F'(L)-4f(t))P_{ij}+(2f^2(t)-F'(L)f(t)+f'(t))g_{ij}
\end{split}
\end{equation}
In the last step we have used the condition that the sectional curvature is nonnegative and the Ricci curvature is parallel. Now in view of the conditions that $F$ is convex and $$2f^2(t)-F'(L)f(t)+f'(t)\geq 0,$$
we have that
$$
\frac{\partial}{\partial t}P_{ij}
\geq\Delta P_{ij}+2\nabla_lL\cdot\nabla_lP_{ij}+2R_{kijl}P_{kl}-R_{il}P_{jl}-R_{jl}P_{il}+2P_{il}P_{jl}+(F'(L)-4f(t))P_{ij}.
$$
By using maximum principle of Hamilton for tensor, we complete the proof of the Theorem 1.5.

 {\bf Proof of Theorem 1.6} Set $$P_{i{\bar j}}=\nabla_i\nabla_{\bar j} L+f(t)g_{i{\bar j}}.$$
  It follows from Lemma 3.1 that
  \begin{equation}
\begin{split}
\frac{\partial}{\partial t}P_{ij}=&\frac{\partial}{\partial t}\nabla_i\nabla_{\bar j} L+f'(t)g_{i{\bar j}}\\
=&\Delta\left(\nabla_i\nabla_{\bar{j}}L\right)+R_{i\bar{j}l\bar{k}}\nabla_k\nabla_{\bar{l}}L
+R_{i\bar{j}k\bar{l}}{\nabla_{l}L\nabla_{\bar{k}}L}+\nabla_{k}L\nabla_{\bar{k}}\nabla_i\nabla_{\bar{j}}L\\
&+\nabla_{\bar{k}}L\nabla_k\nabla_i\nabla_{\bar{j}}L
+\nabla_{\bar{j}}\nabla_{k}L\nabla_i\nabla_{\bar{k}}L+\nabla_i\nabla_{k}L\nabla_{\bar{j}}\nabla_{\bar{k}}L\\
&-\frac{1}{2}\left(R_{l\bar{j}}\nabla_i\nabla_{\bar{l}}L+R_{i\bar{l}}\nabla_l\nabla_{\bar{j}}L\right)+F'(L)\nabla_i\nabla_{\bar{j}} L+F''(L)\nabla_iL\nabla_{\bar{j}} L+f'(t)g_{i\bar{j}}\\
=&\Delta P_{i\bar{j}}+R_{i\bar{j}l\bar{k}}P_{k\bar{l}}
+R_{i\bar{j}k\bar{l}}{\nabla_{l}L\nabla_{\bar{k}}L}+\nabla_{k}L\nabla_{\bar{k}}P_{i\bar{j}}\\
&+\nabla_{\bar{k}}L\nabla_kP_{i\bar{j}}
+\left(P_{i\bar{k}}-f(t)g_{i\bar{k}}\right)\left(P_{k\bar{j}}-f(t)g_{k\bar{j}}\right)+\nabla_i\nabla_{k}L\nabla_{\bar{j}}\nabla_{\bar{k}}L\\
&-\frac{1}{2}\left(R_{l\bar{j}}P_{i\bar{l}}+R_{i\bar{l}}P_{l\bar{j}}\right)+F'(L)P_{i\bar{j}}-F'(L)f(t)g_{i\bar{j}}+F''(L)\nabla_iL\nabla_{\bar{j}} L+f'(t)g_{i\bar{j}}
\end{split}
\end{equation}
Since $F$ is a convex function and the holomorphic bisectional curvature is nonnegative, we have
\begin{equation}
\begin{split}
\frac{\partial}{\partial t}P_{ij}\geq
&\Delta P_{i\bar{j}}
+\nabla_{k}L\nabla_{\bar{k}}P_{i\bar{j}}
+\nabla_{\bar{k}}L\nabla_kP_{i\bar{j}}+R_{i\bar{j}l\bar{k}}P_{k\bar{l}}+P_{i\bar{k}}P_{k\bar{j}}+(F'(L)-2f(t))P_{i\bar{j}}\\
&-\frac{1}{2}\left(R_{l\bar{j}}P_{i\bar{l}}+R_{i\bar{l}}P_{l\bar{j}}\right)+\left(f^2(t)-F'(L)f(t)+f'(t)\right)g_{i\bar{j}}
\end{split}
\end{equation}
With the condition (\ref{fconditionKahler}) in hand, we have that
\begin{equation}
\begin{split}
\frac{\partial}{\partial t}P_{ij}\geq
&\Delta P_{i\bar{j}}
+\nabla_{k}L\nabla_{\bar{k}}P_{i\bar{j}}
+\nabla_{\bar{k}}L\nabla_kP_{i\bar{j}}+R_{i\bar{j}l\bar{k}}P_{k\bar{l}}+P_{i\bar{k}}P_{k\bar{j}}\\
&+(F'(L)-2f(t))P_{i\bar{j}}-\frac{1}{2}\left(R_{l\bar{j}}P_{i\bar{l}}+R_{i\bar{l}}P_{l\bar{j}}\right).
\end{split}
\end{equation}
According to the tensor maximum principle of Hamilton, we complete the proof of Theorem 1.6.


\begin{thebibliography}{}


\bibitem{A}Andrews, B.: Harnack inequalities for evolving hypersurfaces. Math. Z. {\bf 217}, 179-197 (1994)

\bibitem{B}Bando, S.: On classification of three-dimensional compact K\"{a}hler manifolds of nonnegative bisectional curvature. J. Diff. Geom. {\bf 19}, 283-297 (1984)

\bibitem{C1}Cao, H. D.:  The K\"{a}hler-Ricci flow on Fano manifolds.
Lecture Notes in Mathematics {\bf 2086}, 239-297 (2013)

\bibitem{C2}Cao, H. D.:  On Harnack's inequalities for the K\"{a}hler-Ricci flow.
Invent. Math. {\bf 109}, 247-263 (1992)

\bibitem{CN}Cao, H. D., Ni, L.: Matrix Li-Yau-Hamilton estimates for the heat equation on K\"{a}hler manifolds.
 Math. Ann. {\bf 331}, 795-807 (2005)

\bibitem{CFL}Cao, X. D., Fayyazuddin Ljungberg B., Liu, B. W.: Differential Harnack estimates for a nonlinear heat equation. J. Funct. Anal. {\bf 265}, 2312-2330 (2013)

\bibitem{Ch1}Chow, B.: On Harnack inequality and entropy for the Gaussian curvature flow. Comm. Pure Appl. Math. {\bf 44}, 469-483 (1991)

\bibitem{Ch3}Chow, B.: The Yamabe flow on locally conformal flat manifolds with positive Ricci curvature. Comm. Pure Appl. Math. {\bf 45}, 1003-1014 (1992)

\bibitem{Ch4}Chow, B.: Interpolating between Li-Yau's and Hamilton's Harnack inequalities on surface. J. Part. Diff. Equ. {\bf 11}, 137-140 (1998)

\bibitem{CC1}Chow, B., Chu, S.C.: A geometric interpretation of Hamilton's Harnack inequality for the Ricci flow. Math. Res. Lett. {\bf 2}, 701-718 (1995)

\bibitem{CC2}Chow, B., Chu, S.C.: A geometric approach to the linear trace Harnack inequality for the Ricci flow. Math. Res. Lett. {\bf 3}, 549-568 (1996)

\bibitem{CH}Chow, B., Hamilton, R.: Constrained and linear Harnack inequalities for parabolic equations. Invent. Math. {\bf 129}, 213-238 (1997)

\bibitem{CK}Chow, B., Knopf, D.: New Li-Yau-Hamilton inequalities for the Ricci flow via the space-time approach. J. Diff. Geom. {\bf 60}, 1-54 (2002)

\bibitem{CLN}Chow, B., Lu, L. and Ni L.: Hamilton's Ricci Flow, Science Press and American Mathematical Society, 2006.

\bibitem{H1}Hamilton, R.: Three-manifolds with positive Ricci curvature.  J. Diff. Geom. {\bf 17}, 255-306 (1982)

\bibitem{H2}Hamilton, R.: Four-manifolds with positive curvature operator.  J. Diff. Geom. {\bf 24}, 153-179 (1986)

\bibitem{H3}Hamilton, R.: The Ricci flow on surfaces. Contemp. Math. {\bf 71}, 237-262 (1993)

\bibitem{H4}Hamilton, R.: A matrix Harnack estimate for the heat equation. Comm. Anal. Geom. {\bf 1}, 113-126 (1993)

\bibitem{H5}Hamilton, R.: The Harnack estimate for the Ricci flow. J. Diff. Geom. {\bf 37}, 225-243 (1993)

\bibitem{H6}Hamilton, R.: Monotonicity formulas for parabolic flows on manifolds. Comm. Anal. Geom. {\bf 1}, 127-137 (1993)

\bibitem{H7}Hamilton, R.: The Harnack estimate for the mean curvature flow. J. Diff. Geom. {\bf 41}, 215-226 (1995)

\bibitem{K}Kotschwar, B: Some results on the qualitative behavior of solutions to the Ricci flow and other geometric evolution equations.
Thesis (Ph.D.) University of California, San Diego. 2007.


\bibitem{LY}Li, P., Yau, S.-T.: On the parabolic kernel of the Schr\"{o}dinger operator. Acta Math. {\bf 156}, 153-201 (1986)


\bibitem{M2}Mok, N.: The uniformization theorem for compact K\"{a}hler manifolds of nonnegative holomorphic bisectional curvature. J. Diff. Geom. {\bf 27}, 179-214 (1988)

\bibitem{N}Ni, L.: A matrix Li-Yau-Hamilton estimate for K\"{a}hler-Ricci flow. J. Diff. Geom. {\bf 75}, 303-358 (2007)


\bibitem{N2}Ni, L.: Monotonicity and K\"{a}hler-Ricci flow. Contemp. Math. {\bf 367}, 149-165 (2005)


\bibitem{N4}Ni, L.: A monotonicity formula on complete K\"{a}hler manifold with nonnegative bisectional curvature. J. Amer. Math. Soc. {\bf 17}, 909-946 (2004)

\bibitem{NT1}Ni, L., Tam, L.-F.: Plurisubharmonic functions and the K\"{a}hler-Ricci flow. Amer. J. Math. {\bf 125}, 623-654 (2003)

\bibitem{NT2}Ni, L., Tam, L.-F.: Plurisubharmonic functions and the structure of complete K\"{a}hler manifolds with nonnegative curvature. J. Diff. Geom. {\bf 64}, 457-524 (2003)

\bibitem{RYSZ}Ren, X.-A., Yao, S., Shen, L.-J., Zhang, G.-Y.: Constrained matrix Li-Yau-Hamilton estimates on K\"{a}hler manifold. Math. Ann. {\bf 361}, 927-941 (2015)

\bibitem{W}Wu, J. Y., New differential Harnack inequalities for nonlinear heat equations. arXiv:1803.10622v1.

\end{thebibliography}


\end{document}